\documentclass[12pt,draft]{amsart}

\usepackage{amsmath}
\usepackage{amssymb}
\usepackage{enumerate}
\usepackage{euscript}
\usepackage{amscd}
\usepackage[all]{xy}
\usepackage{eq}

\newtheorem{thm}{Theorem}[section]
\newtheorem{lem}[thm]{Lemma}
\newtheorem{prop}[thm]{Proposition}
\newtheorem{cor}[thm]{Corollary}
\theoremstyle{definition}
\newtheorem{defn}[thm]{Definition}
\newtheorem{rem}[thm]{Remark}
\newtheorem{conj}[thm]{Conjecture}

\theoremstyle{remark}

\newcommand{\bba}{{\mathbb A}}

\newcommand{\bbc}{{\mathbb C}}

\newcommand{\bbh}{{\mathbb H}}

\newcommand{\bbq}{{\mathbb Q}}
\newcommand{\bbr}{{\mathbb R}}

\newcommand{\bbz}{{\mathbb Z}}

\newcommand{\Gam}{{\Gamma}}
\newcommand{\del}{{\delta}}

\newcommand{\Del}{{\Delta}}

\newcommand{\lam}{{\lambda}}

\newcommand{\lamb}{\underline{\lambda}}

\newcommand{\gC}{{\mathfrak C}}

\newcommand{\gii}{{\mathfrak e}}
\newcommand{\gE}{{\mathfrak E}}
\newcommand{\gf}{{\mathfrak f}}
\newcommand{\gF}{{\mathfrak F}}

\newcommand{\gM}{{\mathfrak M}}

\newcommand{\gp}{{\mathfrak p}}

\newcommand{\gr}{{\mathfrak r}}

\newcommand{\esN}{{\EuScript{N}}}

\newcommand{\cA}{{\mathcal A}}

\newcommand{\cK}{{\mathcal K}}

\newcommand{\cS}{{\mathcal S}}
\newcommand{\cR}{{\mathcal R}}
\newcommand{\co}{{\mathcal O}}
\newcommand{\cD}{{\mathcal D}}

\newcommand{\sgn}{{\operatorname {sgn}}}

\newcommand{\n}{{\operatorname {N}}}

\newcommand{\re}{{\operatorname {Re}}}
\newcommand{\gal}{{\operatorname{Gal}}}
\newcommand{\gl}{{\operatorname{GL}}}

\newcommand{\Res}{{\operatorname{Res}}}

\newcommand{\op}{{\operatorname{op}}}

\newcommand{\sst}{{\operatorname{ss}}}

\newcommand{\vol}{{\operatorname{vol}}}
\newcommand{\ord}{{\operatorname{ord}}}
\newcommand{\sym}{{\operatorname{Sym}}}

\newcommand{\bs}{Schwartz--Bruhat function}

\newcommand{\supp}{\operatorname{supp}}

\newcommand{\A}{\bba}
\newcommand{\Z}{\bbz}
\newcommand{\Q}{\bbq}
\newcommand{\R}{\bbr}
\newcommand{\C}{\bbc}

\newcommand{\ma}{\bba^{\times}}
\newcommand{\mr}{\R_+}

\newcommand{\mk}{k^{\times}}
\newcommand{\ag}{G_{\bba}}

\newcommand{\md}{d^{\times}\!}
\newcommand{\gMf}{\gM_{\rm f}}
\newcommand{\gMi}{\gM_{\infty}}
\newcommand{\gMc}{\gM_\C}
\newcommand{\gMr}{\gM_\R}
\newcommand{\gMdy}{\gM_{\rm {dy}}}
\newcommand{\gMrm}{\gM_{\rm {rm}}}
\newcommand{\gMin}{\gM_{\rm {in}}}
\newcommand{\gMsp}{\gM_{\rm {sp}}}

\newcommand{\ti}{\widetilde}

\newcommand{\cleq}{\preccurlyeq}

\newcommand{\nr}{\EuScript N}
\newcommand{\trc}{\EuScript T}
\newcommand{\Cor}{{\operatorname{Cor}}}
\newcommand{\prj}{\gr}
\newcommand{\twtw}[4]{\left(\begin{array}{cc}{#1}&{#2}\\{#3}&{#4}\\\end{array}\right)}
\newcommand{\lmb}[1]{\underline{\lambda_{#1}}}
\newcommand{\esQ}{{\EuScript{Q}}}
\newcommand{\esB}{{\EuScript{B}}}
\begin{document}

\title[square of class numbers]
{A mean value theorem for the square of class number
times regulator of quadratic extensions}
\author[Takashi Taniguchi]{Takashi Taniguchi}
\address{Graduate School of Mathematical Sciences\\ University of Tokyo\\
3--8--1 Komaba Meguro-ku\\ Tokyo 153-0041\\ JAPAN}
\email{tani@ms.u-tokyo.ac.jp}
\date{\today}
\begin{abstract}
Let $k$ be a number field.
In this paper, we give a formula for the mean value of the square of
class number times regulator for certain families of 
quadratic extensions of $k$ characterized by finitely many local conditions.
We approach this by using the theory of the zeta function 
associated with the space of pairs of quaternion algebras.
We also prove an asymptotic formula of the correlation coefficient
for class number times regulator of certain families of quadratic extensions.
\end{abstract}

\maketitle

\section{Introduction}\label{sec:int}

We fix an algebraic number field $k$.
Let $\gM$,
$\gM_{\infty}$ and $\gM_{\text{f}}$
denote respectively the set of all places of $k$, all infinite
places and all finite places.
For $v\in\gM$ let $k_v$
denotes the completion of $k$ at $v$ and
if $v\in\gMf$ then let $q_v$ denote
the order of the residue field of $k_v$.
We let $\Delta_k$, $r_1$, $r_2$, and $e_k$ be
respectively the absolute discriminant,
the number of real places, the number of complex places,
and the number of roots of unity contained in $k$.
We denote by $\zeta_k(s)$ the Dedekind zeta function of $k$.

Let $S\supset \gM_\infty$ be a finite set of places.
We fix an $S$-tuple $L_S=(L_v)_{v\in S}$ where
each $L_v$ is a separable quadratic algebra of $k_v$,
i.e., either $k_v\times k_v$ or a quadratic extension of $k_v$.
Let $\esQ(L_S)$ be the following family of quadratic extensions of $k$;
\begin{equation*}
\esQ(L_S)
:=\{F\mid [F:k]=2, F\otimes k_v\cong L_v\ \text{for all}\ v\in S\}.
\end{equation*}
Let $h_F$ and $R_F$ be the class number and the regulator of $F$, respectively.
We would like to understand the value
$h_F^2R_F^2$ for $F\in\esQ(L_S)$ in average.
For $F\in\esQ(L_S)$
we denote by $\Delta_{F/k}$ the relative discriminant of $F/k$
and by $\esN(\Delta_{F/k})$ its absolute norm.
Let
\[
\esQ(L_S,X):=\{F\in\esQ(L_S)\mid \esN(\Delta_{F/k})\leq X\}.
\]
The following is one of the main results of this paper.
\begin{thm}[Theorem \ref{maintheorem}]
\label{introtheorem}
Let $L_S=(L_v)_{v\in S}$ be an $S$-tuple such that
$L_v$ is a field for at least two places of $S$.
Then the limit
\begin{equation*}
\lim_{X\to\infty}\frac{1}{X^2}
\sum_{F\in\esQ(L_S,X)}h_F^2R_F^2
\end{equation*}
exists, and the value is equal to
\begin{equation*}
	\frac{(\Res_{s=1}\zeta_k(s))^3\Delta_k^2e_k^2\zeta_k(2)^2}
		{2^{r_1+r_2+1}2^{2r_1(L_S)}(2\pi)^{2r_2(L_S)}}
		\prod_{v\in S\cap\gMf}\gii_v(L_v)
		{\prod_{v\in\gMf}}(1-3q_v^{-3}+2q_v^{-4}+q_v^{-5}-q_v^{-6}).
\end{equation*}
\end{thm}
Here we denote by $r_1(L_S)$ and $r_2(L_S)$ be respectively
the number of real and complex places of $F\in\esQ(L_S,X)$
(which does not depend on the choice of $F$,)
and also for $v\in\gMf$ we put
\begin{equation*}
\gii_v(L_v)
=
\begin{cases}
	2^{-1}(1+q_v^{-1})(1-q_v^{-2})
		&	L_v\cong k_v\times k_v,\\
	2^{-1}(1-q_v^{-1})^3
		&	L_v\text{ is quadratic unramified},\\
	2^{-1}\esN(\Delta_{L_v/k_v})^{-1}(1-q_v^{-1})(1-q_v^{-2})^2
		&	L_v\text{ is quadratic ramified}.\\
\end{cases}
\end{equation*}
We discuss on the condition of $L_S$ in Remark \ref{unconditional}.

We note that there is a good deal of works on moments of
$h_FR_F$ of quadratic fields $F$ over $\Q$.
For example, Granville and Soundararajan \cite{grso} 
recently obtained the mean value of a general complex
power of $h_FR_F$ of quadratic fields $F$.
On the other side little is known explicitly
over a general number field $k$ to the present
except for the result of Datskovsky \cite{dats} and Kable-Yukie \cite{kayuI}.
Datskovsky \cite{dats} has obtained
the mean value of $h_FR_F$ of quadratic extensions.

We explain one more theorem we prove in this paper.
As we determine every constants in Theorem \ref{introtheorem} explicitly,
combined with the result of Kable-Yukie \cite{kayuI},
we can obtain an interesting formula of
the asymptotic behavior of the correlation coefficients
for class number times regulator of certain families of quadratic extensions.

We fix a quadratic extension $\ti k$ of $k$.
Let $\gM_{\text{rm}}$, $\gM_{\text{in}}$ and $\gM_{\text{sp}}$
be the sets of finite places of $k$ which are
respectively ramified, inert and split on extension to $\ti k$.
We assume $\gM_{\text{rm}}$ does not contain places those dividing $2$.
For any quadratic extension $F$ of $k$ other than $\ti k$,
the compositum $F$ and $\ti k$ contains exactly
three quadratic extensions of $k$.
Let $F^\ast$ denote the quadratic extension other than $F$ and $\ti k$.
Take any $F\in\esQ(L_S)$ and put $L_v^\ast=F^\ast\otimes k_v$,
which does not depend on the choice of $F$.
\begin{thm}[Theorem \ref{cor}]\label{introcor}
Assume $S\supset \gMrm\cup\gM_\infty$.
Let $L_S=(L_v)_{v\in S}$ be an S-tuple.
Assume two of $L_v$'s  and two of $L_v^\ast$'s are fields.
Then the limit
\begin{equation*}
\lim_{X\to\infty}
	\frac{\sum_{F\in\esQ(L_S,X)}h_FR_Fh_{F^\ast}R_{F^\ast}}
	{\left(\sum_{F\in\esQ(L_S,X)}h_F^2R_F^2\right)^{1/2}
	\left(\sum_{F\in\esQ(L_S,X)}h_{F^\ast}^2 R_{F^\ast}^2\right)^{1/2}}
\end{equation*}
exists, and the value is equal to
\begin{equation*}
\prod_{v\in\gMin\setminus S}
	\left(1-\frac{2q_v^{-2}}{1+q_v^{-1}+q_v^{-2}-2q_v^{-3}+q_v^{-5}}\right).
\end{equation*}
\end{thm}
It is an interesting phenomenon that the value
is purely of product form with the index set $\gMin$.
For example, if we take $\ti k$ such that 
$\ti k$ splits at all the small places of $k$,
then $h_FR_F$ and $h_{F^\ast}R_{F^\ast}$ are strongly related.

\smallskip

We prove these density theorems using the theory of {\em zeta functions}
associated with {\em prehomogeneous vector spaces}.
This method has been developed by Sato-Shintani \cite{sash},
Shintani \cite{shintanib},
Datskovsky-Wright \cite{dawrb}, Datskovsky \cite{dats}
and also by Kable, Yukie and the author.
In the beautiful work of Wright-Yukie \cite{wryu},
they showed that 8 types of prehomogeneous vector space
possess significant interest in arithmetic,
and laid out a program to prove a series of density theorems.
There are some advantages to using this theory.
For example, at the moment this approach is
the only possible way that allows the ground field
to be a general number field rather than just $\Q$,
as is done in \cite{dawrb}, \cite{dats, pcbq}, \cite{kayuI} or \cite{ddca}.

This paper is concerned with the representation
\begin{equation*}
G'=\gl(2)\times\gl(2)\times\gl(2),\quad
V'=k^2\otimes k^2\otimes k^2,
\end{equation*}
which is referred to as the $D_4$ case in \cite{wryu}.
It was found in \cite{wryu} that the principal parts of the zeta function
of this type are closely related to
the asymptotic behavior of the mean value of $h_F^2R_F^2$
of quadratic extensions $F/k$.
However, the global theory of prehomogeneous vector spaces
is difficult in general
and more than ten meaningful cases
including the case $(G',V')$ are left open.

Our approach to work on this topic is to consider inner forms.
Let $\esB$ be a quaternion algebra of $k$ and
$\esB^{\rm op}$ the opposite algebra.
We regard $\esB^\times$ and $(\esB^{\rm op})^\times$
as algebraic groups over $k$.
In this paper, we consider the representation
\begin{equation*}
G=\esB^\times\times(\esB^{\rm op})^\times\times\gl(2),\quad
V=\esB\otimes k^2,
\end{equation*}
which is an inner form of $(G,V)$.
Note that if $\esB$ splits then $(G,V)$ is equivalent to $(G',V')$.
We call $(G,V)$ the {\em space of pairs of quaternion algebras}.
As we saw in \cite{zfps},
the orbit space of $V$ also carries a rich structure.
We recall the fundamental properties
of this space in Section \ref{sec:rev}.
One advantage of non-split cases is
that the global theory becomes much easier.
In this paper we consider $(G,V)$
when $\esB$ is a division algebra over a number field $k$.
For this case, we determined the principal parts
of the global zeta function in \cite{zfps}. 

On the other hand, as we will see in Proposition \ref{ZetaSeries},
the global zeta function is only an approximation
of the counting function of $h_F^2R_F^2$ of quadratic extensions.
Hence we could not directly deduce Theorem \ref{introtheorem}
from the global theory \cite{zfps},
and what the aim of this paper is to fill out this gap
by carrying out what is called the {\em filtering process}
originally developed by Datskovsky-Wright \cite{dawra,dawrb}
and Datskovsky \cite{dats}.
This process requires a local theory in some detail.
We consider the localizations of $(G,V)$ at each place of $k$.
Except for a finite set of places of $k$ the quaternion $\esB$ splits,
and the localizations of $(G,V)$ 
at those places are equivalent to $(G',V')$.

There also exists an outer form of the representation $(G',V')$,
namely {\em the space of pairs of binary Hermitian forms}.
The necessary local theory and the filtering process
for that case was constructed
by Kable-Yukie \cite{kayuI, kayuII, kayuIII},
and certain new density theorems were obtained
using the Yukie's global theory \cite{yukieh}.
Some results obtained in \cite{kayuI, kayuII} are useful for us,
because they also consider the split form $(G',V')$ in local situations.
We quote from \cite{kayuI, kayuII}
a uniform estimate of the standard local zeta functions
and some evaluated constants.

After we prove Theorem \ref{maintheorem},
we consider on the correlation coefficient in the final section
combined with the results of \cite{kayuI}.
It is an interesting phenomenon that
the density is purely of the Euler product form
as is stated in Theorem \ref{introcor}.

In \cite{zfps} we handled one more prehomogeneous vector space
which is a non-split form of so called $E_6$-type.
The density theorem for that case,
should be the mean value of $h_FR_F$ of certain families of cubic 
extensions $F$ of $k$, will be treated in a separate paper.
Also in \cite{zfps} we develop the global theory with
general quasi-characters but in this paper we only treat
the principal quasi-character case.
The study with non-principal characters will enrich
arithmetic results as is done in \cite{dawra}, \cite{kawr}, or \cite{ddca}.
We hope this to be developed in the future.
\smallskip

For the remainder of this section, we will give the contents of the paper.
In Section \ref{sec:not}, we introduce the notations used throughout the paper.
More specialized notations are introduced when required.
In Section \ref{sec:rev}, we define the space of pairs of quaternion algebras,
and recall from \cite{zfps} its basic properties.
In Section \ref{sec:inv}, we first define various invariant measures
on the groups and the representation spaces.
After that we introduce the global zeta function
and review its analytic properties.

From Section \ref{sec:can} to Section \ref{sec:inf},
we consider the local theory.
We establish the necessary local theory to obtain the density theorem
in these sections.
In Section \ref{sec:can}, we define a measure on the stabilizer
for semi-stable points, which is in some sense canonical.
In Section \ref{sec:loc}, we define the local zeta function
and the local density.
Also we quote from \cite{kayuI}
an estimate of the standard local zeta function,
which we need in order to apply the filtering process
in the proof of the mean value theorem in Section \ref{sec:mea}.
In Sections \ref{sec:fin}, \ref{sec:ram} and \ref{sec:inf},
we compute the local densities.
Section \ref{sec:fin} is for finite unramified places
(the places $\esB$ splits),
Section \ref{sec:ram} for finite ramified places, and
Section \ref{sec:inf} for infinite places.
The unramified cases were almost done
in \cite{kayuI,kayuII} and we essentially quote their result,
but we will give a refinement for
dyadic places by applying the method developed in \cite{pcbq}. 
After that we study the ramified cases.

In Section \ref{sec:mea} we go back to the adelic situation.
We first define some invariant measures and show that
our zeta function is more or less the counting function
of the unnormalized Tamagawa numbers of the stabilizers.
After that we apply the filtering process to our case
and find the mean value of the Tamagawa numbers.
Then with an explicit computation, we give a formula
for the mean value of the square of
class numbers times regulators for certain family of 
quadratic extensions, which is a main theorem of this paper.
In Section \ref{sec:cor}, we define
the correlation coefficient
of class number times regulator of quadratic fields.
Then we explicitly compute the value in some cases
by combining the results of \cite{kayuI} and this paper.

\bigskip

\noindent
{\bf Acknowledgments.}
There are many tremendous help for the creation of this work.
The author express his gratitude to his advisor T. Terasoma
for the constant discussions and suggestions.
The author also would like to thank to Professor A. Yukie
who suggested to consider the topic in Section \ref{sec:cor},
and to A. C. Kable who taught me several references
including \cite{grso} with useful comments.
The author was also inspired by their
series of work \cite{kayuI,kayuII, kayuIII}.
Special thanks goes to the author's colleague Uuye Otogonbayar,
who read the manuscript and gave many comments.

\section{Notation}\label{sec:not}

In this section we collect basic notations used throughout in this paper.

If $X$ is a finite set then $\# X$ will denote its cardinality. The
standard symbols $\Q$, $\R$, $\C$ and $\Z$ will denote respectively
the rational, real and complex numbers and the rational integers.
The set of positive real numbers is denoted $\R_{+}$.
For a complex number $z$, let $\Re(z),\Im(z)$ and $\bar z$ be
the real part, the imaginary part, and the complex conjugate of $z$.
If $R$ is any ring then $R^{\times}$ is the set
of invertible elements of $R$, and if $V$ is a scheme defined over
$R$ and $S$ is an $R$-algebra
then $V_S$ denotes its $S$-rational points.
Let us denote by ${\rm M}(2,2)$ the set of $2\times 2$ matrices.

We fix an algebraic number field $k$.
Let $\gM$,
$\gM_{\infty}$, $\gM_{\text{f}}$, $\gM_{\text{dy}}$, 
$\gM_{\R}$ and $\gM_{\C}$
denote respectively the set of all places of $k$, all infinite
places, all finite places, all dyadic places (those dividing the
place of $\Q$ at 2), all real places and all complex places.
Let $\co$ be the ring of integers of $k$. If $v\in\gM$ then $k_v$
denotes the completion of $k$ at $v$ and $|\;|_v$ or $|\;|_{k_v}$ 
denotes the normalized
absolute value on $k_v$. If $v\in\gMf$ then $\co_v$ denotes
the ring of integers of $k_v$,
$\gp_v$ the maximal ideal of $\co_v$ and $q_v$ the cardinality of
$\co_v/\gp_v$. 
For $t\in\ k_v^\times$, we define $\ord_v(t)$ so that
$|t|_v=q_v^{-\ord_v(t)}$.
For a practical purpose in Sections \ref{sec:fin} and \ref{sec:ram},
we do {\em not} fix a uniformizer in $\co_v$ here.
For any separable quadratic algebra $L_v$ of $k_v$,
let $\co_{L_v}$ denote the ring of integral elements of $L_v$.
That is, if $L_v$ is a quadratic extension then $\co_{L_v}$ is
the integer ring of $L_v$ and 
if $L_v=k_v\times k_v$ then $\co_{L_v}=\co_v\times\co_v$.

If $k_1/k_2$ is a finite extension of either local fields or
number fields then we shall write $\Del_{k_1/k_2}$ for the relative
discriminant of the extension; it is an ideal in the ring of
integers of $k_2$.
For conventions, we let $\Del_{k_2\times k_2/k_2}$ be
the integer ring of $k_2$.
If the extension $k_1/k_2$ is of number fields,
let $\esN(\Delta_{k_1/k_2})$ be the absolute norm of $\Delta_{k_1/k_2}$.
The symbol $\Del_{k_1}$ will stand for $\esN(\Delta_{k_1/\Q})$,
the classical absolute discriminant of $k_1$ over $\Q$.
We use the notation $\n_{k_1/k_2}$ for
the norm in $k_1/k_2$.

Returning to $k$, we let $r_1$, $r_2$, $h_k$, $R_k$ and $e_k$ be
respectively the number of real places, the number of complex
places, the class number, the regulator and the number of roots of
unity contained in $k$. It will be convenient to set
\begin{equation*}
\gC_k=2^{r_1}(2\pi)^{r_2}h_kR_ke_k^{-1}
.\end{equation*}

We refer to \cite{weilc} as the basic reference for fundamental properties
on adeles.
The ring of adeles and the group of ideles
are denoted by $\A$ and $\ma$, respectively.
The adelic absolute value $|\;|$ on $\ma$ is normalized so that,
for $t\in\ma$, $|t|$ is the module of multiplication by $t$
with respect to any Haar measure $dx$ on $\A$, i.e. $|t|=d(tx)/dx$.
Let $\A^0=\{t\in\ma\mid |t|=1\}$.
Suppose $[k:\Q]=n$.
For $\lam\in\R_+$, $\lamb\in \ma$ 
is the idele whose component at any infinite place is $\lam^{1/n}$
and whose component at any finite place is $1$.  
Then we have $|\lamb|=\lam$.

For a finite extension $L/k$,
let $\A_L$ denote the adele ring of $L$.
We define $\A_L^\times,\A_L^0,\gC_L$ etc., similarly.
The adelic absolute value of $L$ is denoted by $|\;|_L$.
There is a natural inclusion $\A\to \A_L$, 
under which an adele $(a_v)_v$ corresponds to the adele
$(b_w)_w$ with $b_w=a_v$ if $w|v$. 
Using the identification $L\otimes_{k} \A \cong \A_L$,
the norm map 
$\n_{L/k}$ can be extended to a map from $\A_L$ to $\A$.  
It is known (see p. 139 in \cite{weilc}) that 
$|\n_{L/k}(t) | = |t|_{L}$ for $t\in \ti \A$. 
Suppose $[L:k]=m$.
For $\lam\in\mr$,
we denote by $\lamb_L\in \A_L^{\times}$ 
are the ideles whose component at 
any infinite place is $\lam^{1/mn}$
and whose component at any finite place is $1$,
so that $|\lamb_L|_L=\lam$.
Clearly $\lamb = \lamb_L^m$
and hence $|\lamb|_L= \lam^m$.  When we have to show the 
number field on which we consider $\lamb$, we use the notation 
such as $\lamb_k$.

If $V$ is a vector space over $k$ we write $V_{\A}$ for its adelization.
Let $\cS(V_{\A})$ and $\cS(V_{k_v})$
be the spaces of \bs s on each of the indicated domains.

For any  $v\in \gM_{\text{f}}$, 
we choose a Haar measure $dx_v$ on $k_v$ 
to satisfy $\int_{\co_v}dx_v=1$.  
We write $dx_v$ for the ordinary Lebesgue measure if $v$ is real, 
and for twice the Lebesgue measure if $v$ is imaginary.  
We choose a Haar measure  $dx$ on $\A$ to satisfy $dx =\prod_{v\in\gM} dx_v$.
Then $\int_{\A/k}dx=|\Delta_k|^{1/2}$
(see \cite{weilc}, p. 91).  

For any $v\in \gM_{\text{f}}$, we normalize the Haar measure $\md t_v$
on $\mk_v$ such that $\int_{\co_v^{\times}} \md t_v = 1$.  
Let $\md t_v(x)=|x|_v^{-1}d x_v$ if $v\in\gMi$.
We choose a Haar
measure $\md t$ on $\ma$  so that $\md t=\prod_{v\in\gM}\md t_v$.
Using 
this measure,
we choose a Haar measure $\md t^0$ on $\A^0$ by
$$
\int_{\ma} f(t) \,\md t  = \int_0^{\infty}   
\int_{\A^0}  f(\lamb t^0) \,\md \lam \md t^0
,$$ 
where $\md \lam = \lam^{-1}d\lam$.
Then $\int_{\A^0/k^\times}\md t^0=\gC_k$  
(see \cite{weilc}, p. 95). 

Let $\zeta_k(s)$ be the Dedekind zeta function of $k$.
We define
\begin{equation*}\label{Zkdefn}
Z_k(s)=|\Delta_k|^{s/2}
\left(\pi^{-s/2}\Gam\left(\frac s 2\right)\right)^{r_1}
\left((2\pi)^{1-s}\Gam(s)\right)^{r_2}\zeta_k(s)\,
.\end{equation*}
This definition differs from that in \cite{weilc}, p.129 by the
inclusion of the $|\Del_{k}|^{s/2}$ factor.
It is adopted here as
the most convenient for our purposes.
It is known 
(\cite{weilc}, p.129) that 
\begin{equation*} \label{zetaresidues} 
\operatorname{Res}_{s=1} \zeta_k(s) = 
|\Delta_k|^{-{\frac 12}}\gC_k\text{\quad and so\quad}
\operatorname{Res}_{s=1} Z_k(s) = \gC_k\,
.\end{equation*}

Let $\bbh$ denote the quaternion algebra of Hamiltonians over $\R$.
We choose and fix an element $j\in\bbh$ so that
$\bbh=\C\oplus\C j$ as a left vector space over $\C$
and the multiplication law is given by
$j^2=-1$ and 
$j\alpha=\bar\alpha j$ for $\alpha\in\C$.
Let us express elements of $\bbh$ as $x=x_1+x_2j$ where $x_1,x_2\in\C$.
We choose a Haar measure on $\bbh$
so that $dx=dx_1dx_2$,
where $dx_1$ and $dx_2$ are twice the Lebesgue measure on $\C$
as above. If we let $|x|_\bbh=|x_1|_\C+|x_2|_\C$,
then $|x|_\bbh^{-2}dx$ defines a Haar measure on $\bbh^\times$.
For practical purposes, we choose $\md t(x)=\pi^{-1}|x|_\bbh^{-2}dx$
as the normalized measure on $\bbh^\times$.
We put $\bbh^0=\{t\in\bbh^\times\mid |t|_\bbh=1\}$.

\section{Review of the space of pairs of quaternion algebras}\label{sec:rev}

In this section, we define the prehomogeneous vector space
of pairs of quaternion algebras
which are at the heart of this work and reviewing their
fundamental properties.
Arithmetic plays no role here, so in this section
we consider the representation over an arbitrary field $K$.
We later use the result in this section both local and global situations.

Let $\esB$ be a quaternion algebra over $K$.
This algebra is either isomorphic to the algebra
${\rm M}(2,2)$ consisting of $2\times2$ matrices or
a division algebra of dimension $4$.
Let $\trc$ and $\nr$ be the reduced trace and the reduced norm, respectively.
We denote by $\esB^\op$ the opposite algebra of $\esB$.
We introduce a group $G_1$ and its linear representation on $\esB$ as follows.
Let
\begin{equation*}
G_{11}=\esB^\times,\quad G_{12}=(\esB^{\rm op})^\times,
\quad\text{and}\quad
G_1=G_{11}\times G_{12}.
\end{equation*}
That is, $G_1$ is equal to $\esB^\times\times \esB^\times$ set theoretically
and the multiplication law is given by
$(g_{11},g_{12})(h_{11},h_{12})=(g_{11}h_{11},h_{12}g_{12})$.
If there is no confusion, we drop `op' and simply write such as
$G_{12}=\esB^\times$ instead.
We regard $G_1$ as an algebraic group over $K$.
The quaternion algebra $\esB$ can be considered
as a vector space over $K$.
We define the action of $G_1$ on $\esB$ as follows:
\begin{equation*}
(g_1,w)\longmapsto g_{11}w g_{12},
	\qquad g_1=(g_{11},g_{12})\in G_1,
			w\in \esB.
\end{equation*}
This defines a representation $\esB$ of $G_1$.
We consider the standard representation of $G_2=\gl(2)$ on $K^2$.
The group $G=G_1\times G_2$ acts naturally on
$V=\esB\otimes K^2$.
The representation $(G,V)$ is the main object of this paper.
This is a $K$-form of
\begin{equation}\label{splitD4}
(\gl(2)\times\gl(2)\times\gl(2),K^2\otimes K^2\otimes K^2),
\end{equation}
and if $\esB$ is split, $(G,V)$ is equivalent to
the above representation over $K$.
The representation \eqref{splitD4} was studied in \cite{wryu}
in some detail, and our review is
a slight generalization \cite{zfps} of that.

We describe the action more explicitly.
Throughout this paper,
we express elements of $V\cong \esB\oplus \esB$ as $x=(x_1,x_2)$.
We identify $x=(x_1,x_2)\in V$ with $x(v)=v_1x_1+v_2x_2$
which is an element of the quaternion algebra with coordinates
in linear forms in two variables $v=(v_1,v_2)$.
Then the action of $g=(g_{11},g_{12},g_2)\in G$ on $x\in V$ is defined by
\begin{equation*}
(gx)(v)=g_{11}x(vg_2)g_{12}.
\end{equation*}
We put $F_x(v)=\nr(x(v))$.
This is a binary quadratic form in variables $v=(v_1,v_2)$.
We let $P(x)\ (x\in V)$ be the discriminant of $F_x(v)$,
which is a polynomial in $V$.
That is, if we express $F_x(v)=a_0(x)v_1^2+a_1(x)v_1v_2+a_2(x)v_2^2$,
then $P(x)$ is given by $P(x)=a_1(x)^2-4a_0(x)a_2(x)$.
Let $\chi_i\, (i=1,2)$ be the character of $G_i$ defined by
\begin{equation*}
\chi_1(g_1)=\nr(g_{11})\nr(g_{12}),\quad \chi_2(g_2)=\det g_2,
\end{equation*}
respectively.
We define $\chi(g)=\chi_1(g_1)^2\chi_2(g_2)^2$.
Then one can easily see that
\begin{equation*}
P(gx)=\chi(g)P(x)
\end{equation*}
and hence $P(x)$ is a relative invariant polynomial
with respect to the character $\chi$.
Let $V^\sst=\{x\in V\mid P(x)\not=0\}$,
and we call this the set of semi-stable points.
That is, $x\in V$ is semi-stable if and only if
$F_x(v)$ does not have a multiple root in ${\mathbb P}^1=\{(v_1:v_2)\}$.

Let $\ti T=\ker(G\rightarrow \gl(V))$.
Then it is easy to see that
\begin{equation*}
\ti T=\{(t_{11},t_{12},t_2)\mid 
t_{11},t_{12},t_2\in\gl(1), t_{11}t_{12}t_2=1\},
\end{equation*}
which is contained in the center of $G$.
Throughout this paper,
we will identify $\ti T$ with $\gl(1)^2$ via the map
\begin{equation*}
\ti T\longrightarrow \gl(1)^2
\qquad
(t_{11},t_{12},(t_{11}t_{12})^{-1})\longmapsto (t_{11},t_{12}).
\end{equation*}

We are now ready to recall the description of
the space of non-singular $G_K$-orbits in $V_K$.
\begin{defn}
For $x\in V_K^\sst$, we define
\begin{align*}
Z_x		&=\mathrm{Proj}\, K[v_1,v_2]/(F_x(v)),\\
\ti K(x)	&=\Gamma(Z_x,\co_{Z_x}).
\end{align*}
Also we define $K(x)$ to be the splitting field of $F_x(v)$.
\end{defn}
Note that $\ti K(x)$ may not be a field.
Since $V_K^\sst$ is the set of $x$ such that
$F_x$ does not have a multiple root,
$Z_x$ is a reduced scheme over $K$ and
$\ti K(x)$ is a separable commutative $K$-algebra of dimension $2$.
By definition, we immediately see that
$\ti K(x)\cong K\times K$ if $F_x(v)$ has $K$-rational factors and
$\ti K(x)\cong K(x)$ if $F_x(v)$ is irreducible over $K$.

The following lemma is easy to prove,
but quite useful for our practical purposes.
For the proof, see \cite[Lemma 3.3]{zfps}.
\begin{lem}\label{lem:orbit}
Any $G_K$-orbit in $V_K^{\sst}$
contains an element of the form
$w_u=(1,u)$ for some $u\in \esB_K$.
\end{lem}
Note that for $u\in \esB_K$,
$w_u=(1,u)$ is semi-stable if and only if
$u$ is a separable quadratic element of $\esB_K$.
For $w_u\in V_K^\sst$,
$F_{w_u}(v_1,-1)=\nr(v_1-u)$ is the characteristic polynomial of $u$
and hence $\ti K(x)$ is isomorphic to $K[u]\subset \esB_K$ as a $K$-algebra.

\begin{defn}
Let $\cA_2(\esB_K)$ be the set of isomorphism classes of
separable commutative $K$-algebras of dimension $2$
those are embeddable into $\esB_K$.
\end{defn}
Note that if $\esB_K$ is non-split, then any element of $\cA_2(\esB_K)$
is a quadratic extension of $K$.
The following proposition is proved in \cite{wryu}, \cite{zfps}.
\begin{prop}\label{rod}
The map
$x\mapsto \ti K(x)$
gives a bijection between
$G_K\backslash V_K^\sst$ and $\cA_2(\esB_K)$.
\end{prop}
For $x\in V_K^\sst$, let
$G_x$ be the stabilizer of $x$ and
$G_x^\circ$ its identity component,
both are algebraic groups defined over $K$.
We have shown in \cite{zfps} that
$G_x^\circ\cong (\gl(1)_{\ti K(x)})^2$
as an algebraic group over $K$.
We close this section with a detailed description of the
$K$-rational points of the stabilizer $G_{w_u}$.

We first recall the isomorphism $G_{w_u\,K}^\circ\cong (K[u]^\times)^2$.
Since $\{1,u\}$ is a basis of $K[u]$ as a $K$-vector space,
for any $s_1,s_2\in K[u]^\times$,
$\{s_1s_2,s_1s_2u\}$ is also a $K$-basis of $K[u]$.
Hence there exists a unique element
$g_{s_1s_2}\in\gl(2)_K$ such that
$g_{s_1s_2}\,\,^{t}(s_1s_2,s_1s_2u)={}^t(1,u)$.
Since $K[u]$ is a commutative algebra, $s_1s_2u=s_1us_2$.
Therefore we have $(s_1,s_2,g_{s_1s_2})\in G_{w_u\,K}$.
The following proposition is proved in \cite[Lemma 3.4]{zfps}.
\begin{prop}\label{stbstrk}
The map
\begin{equation*}
\psi_u\colon
(K[u]^\times)^2\longrightarrow G_{w_u\,K}^\circ,
\qquad
(s_1,s_2)\longmapsto (s_1,s_2,g_{s_1s_2})
\end{equation*}
gives an isomorphism of the two groups.
\end{prop}
Finally we consider the structure of $G_{w_u\,K}/G_{w_u\,K}^\circ$.
Let $\sigma$ be the non-trivial $K$-automorphism of $K[u]$.
Then there exists $\nu\in \esB_K\setminus K[u]$ such that $\nu^2\in K$,
$\esB_K=K[u]\oplus K[u]\nu$ as a $K[u]$-vector space,
and the multiplication law is given by $\nu\alpha=\alpha^\sigma\nu$
for $\alpha\in K[u]$.
Let $u=a+bu^\sigma$ where $a\in K, b\in K^\times$.
\begin{prop}\label{notidentity}
We have $[G_{w_u\,K}:G_{w_u\,K}^\circ]=2$
and $G_{w_u\,K}/G_{w_u\,K}^\circ$
is generated by the class of
$\tau=\left(\nu^{-1},\nu,\twtw 10ab\right)$.
\end{prop}
\begin{proof}
A simple computation shows $\tau w_u=w_u$.
On the other hand,
by \cite{wryu} we have
$[G_{w_u\,\bar K}:G_{w_u\,\bar K}^\circ]=2$
because $(G,V)$ is a $K$-form of \eqref{splitD4}.
Since $[G_{w_u\,K}:G_{w_u\,K}^\circ]\leq
[G_{w_u\,\bar K}:G_{w_u\,\bar K}^\circ]$,
the proposition follows.
\end{proof}
By Lemma \ref{lem:orbit}, we have
$[G_{x\,K}:G_{x\,K}^\circ]=2$
for any $x\in V_{K}^\sst$.

\section{Invariant measures and the global zeta function}\label{sec:inv}
For the rest of this paper, we assume
$k$ a fixed number field and $\esB$ a non-split quaternion algebra over $k$.
In this section, we define various invariant measures in both
local and adelic situations and summarize the necessary results.
For the proof, see \cite{vign} for example.
In this paper, we always choose the adelic measure
as the product of local measures,
for the conventions of connection
between global and local theory.
After that we introduce the global zeta function of
the prehomogeneous vector space $(G,V)$ 
and recall from \cite{zfps} its most basic analytic properties.

We define $\gM_\esB$ to be the set of places $v$ of $k$
such that $\esB$ is ramified at $v$.
For $v\in\gM$, let $\esB_v$ denote $\esB\otimes_kk_v$.
Then, by definition, $v\in\gM_\esB$ if and only if $\esB_v$ is a division algebra.
It is well known that $\gM_\esB$ is a finite set.

We give a normalization of invariant measure on $G_{k_v}$ and $V_{k_v}$.
First we consider the places $v\notin\gM_\esB$.
For each of these $v$, we fix once and for all a $k_v$-isomorphism
$\esB_v\cong {\rm M}(2,2)_{k_v}$
and identify these algebras.
Then
\begin{equation*}
G_{k_v}=\gl(2)_{k_v}\times\gl(2)_{k_v}\times\gl(2)_{k_v},
\qquad
V_{k_v}={\rm M}(2,2)_{k_v}\oplus{\rm M}(2,2)_{k_v}.
\end{equation*}
We choose a Haar measure $dx_v$ on $V_{k_v}$ so that
\begin{equation*}
dx_v=dx_{1v}dx_{2v},\qquad
dx_{iv}=dx_{i11v}dx_{i12v}dx_{i21v}dx_{i22v}
\ \ (i=1,2)
\end{equation*}
for
\begin{equation*}
x_v=(x_{1v},x_{2v}),\qquad
x_{iv}=\twtw {x_{i11v}}{x_{i12v}}{x_{i21v}}{x_{i22v}}
\ \ (i=1,2).
\end{equation*}
For $v\in\gMf$, we put
$V_{\co_v}={\rm M}(2,2)_{\co_v}\oplus{\rm M}(2,2)_{\co_v}$,
which is a maximal compact subgroup of $V_{k_v}$.
We note that $\int_{V_{\co_v}}dx=1$ for $v\in\gMf$.
We consider $G_{k_v}$ for $v\notin\gM_\esB$.
If $v\in\gMf$, we put a maximal compact subgroup $K_v$ of $G_{k_v}$ as
\begin{equation*}
\cK_v=\gl(2)_{\co_v}\times\gl(2)_{\co_v}\times\gl(2)_{\co_v},
\end{equation*}
and normalize the measure $dg_v$ on $G_{k_v}$ so that
the total volume of $\cK_v$ is $1$.
For $v\in\gMi$,
we first give a measure for $\gl(2)_F$
where $F=\R$ or $\C$.
As in Section \ref{sec:not},
we shall take Lebesgue measure to be the standard measure on the real numbers
and twice the Lebesgue measure to be the standard measure on the complex numbers.
If $h_v=(h_{ijv})_{1\leq i,j\leq2}$,
then $d\mu(h_v)=dh_{11v}dh_{12v}dh_{21v}dh_{22v}/|\det(h_v)|_F^2$
defines a Haar measure on $\gl(2)_F$.
We put $dh_v=p_Fd\mu(h_v)$ where $p_{\R}=\pi^{-1}$ and $p_{\C}=(2\pi)^{-1}$.
Using this measure, we define
$dg_v$ for $v\in\gMi$ as
$dg_v=dg_{11v}dg_{12v}dg_{2v}$ where
$g=(g_{11},g_{12},g_2)\in G_{k_v}=(\gl(2)_{k_v})^3$.

Next we consider the case $v\in\gM_\esB$.
For $v\in\gMf$, let $\co_{\esB_v}$
be the ring consisting of integral elements of $\esB_v$.
We put $V_{\co_v}=\co_{\esB_v}\oplus\co_{\esB_v}$,
which is a maximal compact subgroup of $V_{k_v}$.
We choose a Haar measure $dx_v$ on $V_{k_v}=\esB_v\oplus \esB_v$
so that the volume of $\co_{\esB_v}\oplus\co_{\esB_v}$ is $1$.
Also we we put a maximal compact subgroup $K_v$ of $G_{k_v}$ as
\begin{equation*}
\cK_v=\co_{\esB_v}^\times\times\co_{\esB_v}^\times\times\gl(2)_{\co_v},
\end{equation*}
and normalize the measure $dg_v$ on $G_{k_v}$ so that
the total volume of $\cK_v$ is $1$.

Now the remaining case is for $v\in\gMi\cap\gM_\esB$,
which is an element of $\gM_\R$.
We fix an isomorphism $\esB_v\cong \bbh$. Then
\begin{equation*}
G_{k_v}=\bbh^\times\times\bbh^\times\times\gl(2)_{\R},
\qquad
V_{k_v}=\bbh\oplus\bbh.
\end{equation*}
We set measures $dg_v$ and $dx_v$
on $G_{k_v}$ and $V_{k_v}$ as the product measures,
where we consider the measures on $\bbh^\times,\bbh$
as in Section \ref{sec:not} and $\gl(2)_\R$ as above.
For $v\in\gMi$, we put
\begin{equation*}
\cK_v=
\begin{cases}
{\rm O}(2,\R)^3&v\in\gMr\setminus\gM_\esB,\\
{\rm U}(2,\C)^3&v\in\gMc\setminus\gM_\esB,\\
\bbh^0\times\bbh^0\times{\rm O}(2,\R)&v\in\gM_\esB,\\
\end{cases}
\end{equation*}
which is a maximal compact subgroup of $G_{k_v}$.

Using these local measures, we define the measures
$dg$ and $dx$ on $G_\A$ and $V_\A$ by
\begin{equation*}
dg=\prod_{v\in\gM}dg_v,
\qquad
\text{and}
\qquad
dx=\prod_{v\in\gM}dx_v.
\end{equation*}
If we put
\begin{equation*}
\Delta_\esB=\Delta_k^4\prod_{v\in\gM_\esB\cap\gMf}q_v^2
\qquad
\text{and}
\qquad
\Delta_V=\Delta_\esB^2,
\end{equation*}
then it is well known that the volume of $V_\A/V_k$
with respect to the measure $dx$ is $\Delta_V^{1/2}$.
Hence our choice of measure $dx$ on $V_\A$ in this paper is
$\Delta_V^{1/2}$ times that of \cite{zfps},
in which we defined so that
the volume of $V_\A/V_k$ is equal to $1$.

Our definition of measure $dg_v$ on $G_{k_v}$
can naturally be considered as the product measure
$dg_v=dg_{11v}dg_{12v}dg_{2v}$
for $g_v=(g_{11v},g_{22v},g_{2v})$
and we shall do so below.
For example, if $v\in\gMf\cap\gM_\esB$,
we will regard $dg_{11v},dg_{12v}$ and $dg_{2v}$ 
on $G_{11k_v},G_{12k_v}$ and $G_{2k_v}$ as
\begin{equation*}
\int_{\co_{\esB_v}^\times}dg_{11v}=
\int_{\co_{\esB_v}^\times}dg_{12v}=
\int_{\gl(2)_{\co_v}}dg_{2v}=1.
\end{equation*}
We define the measure $dg_{11},dg_{12}$ and $dg_2$
on $G_{11\A},G_{12\A}$ and $G_{2\A}$ by
\begin{equation*}
dg_{11}=\prod_{v\in\gM}dg_{11v},
\quad
dg_{12}=\prod_{v\in\gM}dg_{12v},
\quad
\text{and}
\quad
dg_{2}=\prod_{v\in\gM}dg_{2v}.
\end{equation*}
Clearly, we have $dg=dg_{11}dg_{12}dg_2$.

Since $\ti T\cong\gl(1)\times\gl(1)$ is a split torus,
the first Galois cohomology set $H^1(k',\ti T)$ is trivial
for any field $k'$ containing $k$.
This implies that the set of $k'$-rational point of $\ti G$
coincides with $G_{k'}/\ti T_{k'}$.
Therefore $(G/T)_\A=G_\A/\ti T_\A$ and
$(G/T)_\A/(G/T)_k=G_\A/\ti T_\A G_k$.
We put the measures $\md \ti t_v$ and $\md\ti t$
on $\ti T_v$ and $\ti T_\A$ respectively to satisfy
$\md \ti t_v=\md t_{1v}\md t_{2v}$, $\md \ti t=\md t_1\md t_2$
for $\ti t_v=(t_{1v},t_{2v},(t_{1v}t_{2v})^{-1}),
\ti t=(t_{1},t_{2},(t_{1}t_{2})^{-1})$.
Using these, we normalize the invariant measure
$d\ti g_v$ and $d\ti g$ on $G_{k_v}/T_{k_v}$ and $G_\A/T_\A$
so that
$dg_v=d\ti g_v\md\ti t_v, dg=d\ti g\md\ti t$.
Note that $d\ti g=\prod_{v\in\gM}d\ti g_v$
since $dg=\prod_{v\in\gM}dg_v$ and $\md\ti t=\prod_{v\in\gM}\md\ti t_v$.

We put
\begin{align*}
G_{1i\A}^0&=\{g_{1i}\in G_{1i\A}\mid |\nr(g_{1i})|=1\}
\ \ (i=1,2),\\
G_{2\A}^0&=\{g_2\in G_{2\A}\mid |\det(g_2)|=1\}.
\end{align*}
Then the maps
\begin{align*}
\R_+\times G_{1i\A}^0\longrightarrow G_{1i\A},
&\quad
(\lam_{1i},g_{1i}^0)\longmapsto \lmb{1i}g_{1i}^0
\ \ (i=1,2),\\
\R_+\times G_{2\A}^0\longrightarrow G_{2\A},
&\quad
(\lam_{2},g_{2}^0)\longmapsto \lmb{2}g_{2}^0,
\end{align*}
give isomorphisms of these groups.
We choose Haar measures $dg_{1i}^0$ and $dg_2^0$
on $G_{1i\A}^0$ and $G_{2\A}^0$
so that
$dg_{1i}=2\md\lam_{1i}dg_{1i}^0,dg_2=2\md\lam_2dg_2^0$.
Then it is known that
\begin{align*}
\int_{G_{1i\A}^0/G_{1ik}}dg_{1i}^0
&=\Delta_k^{1/2}\gC_kZ_k(2)\prod_{v\in\gM_\esB\cap\gMf}(q_v-1),\\
\int_{G_{2\A}^0/G_{2k}}dg_2^0&=\Delta_k^{1/2}\gC_kZ_k(2).
\end{align*}

We now define the global zeta function.

\begin{defn}
For $\Phi\in\cS(V_\A)$ and a complex variable $s$,
we define
\begin{equation*}
Z(\Phi,s)
=	\int_{G_\A/\ti T_\A G_k}
		|\chi(\ti g)|^s\sum_{x\in V_k^\sst}\Phi(\ti gx)\, d\ti g,
\end{equation*}
and call it the {\em global zeta function}.
\end{defn}

It is known that the integral converges if $\Re(s)$ is sufficiently large
and can be continued meromorphically to the whole complex plane.
In \cite{zfps}, we described the principal parts of $Z(\Phi,s)$
by means of certain distributions.
However, we used a slightly different formulation in \cite{zfps},
and we need some arguments to translate the results from that paper.
Also, in this paper we only consider the rightmost pole of $Z(\Phi,s)$
because this is enough to deduce the density theorems.

We put $G_{1\A}^0=G_{11\A}^0\times G_{12\A}^0$.
The domain of integration used in \cite{zfps} is
$\R_+\times G_\A^0/G_k$, where
$G_\A^0=G_{1\A}^0\times G_{2\A}^0$.
Let $\ti T_\A^0=\ag^0\cap \ti T_\A$.
Then we have
\begin{equation*}
(\R_+\times G_\A^0)/\ti T_\A^0\cong G_\A/\ti T_\A
\end{equation*}
via the map which sends the class of
$(\lam,g_{11}^0,g_{12}^0,g_2^0)$ to class of $(g_{11}^0,g_{12}^0,\lamb g_2^0)$.
In \cite{zfps} $\R_+\times G_\A^0$ is made to act on $V_\A$
by assuming that $(\lam,1)$ acts by multiplication by $\lamb$,
and the above isomorphism is compatible with their actions on $V_\A$.
We will compare the measure $d\ti g$ on $G_\A/\ti T_\A$
with the measure $\md\lam dg^0$ on $\R_+\times G_\A^0$ used in \cite{zfps}.
The argument in \cite{zfps} is valid
for any choice of measure on $G_{1\A}^0$
and we consider $dg_{11}^0dg_{12}^0$ for this.
We note that the measure $dg_2^0$ on $G_{2\A}^0$
in the present situation
is $\Delta_k^{1/2}\gC_k^2$ times that of used in \cite{zfps}.

We have
$G_\A/\ti T_\A\cong(\R_+^3\times G_\A^0)/(\R_+^2\times \ti T_\A^0)$
where $\R_+^2\times \ti T_\A^0$ is included in $\R_+^3\times G_\A^0$ via
$(\lam_1,\lam_2,\ti t^0)\mapsto(\lam_1,\lam_2,\lam_1^{-1}\lam_2^{-1},\ti t^0)$
and $\R_+^3\times G_\A^0$ maps onto $G_\A/\ti T_\A$ via
$(\lam_1,\lam_2,\lam_3,g^0)\mapsto(\lmb1,\lmb2,\lmb3)g^0$.
With this identification we have chosen the measure $d\ti g$
to be compatible with the measure
$8\md\lam_1\md\lam_2\md\lam_3 dg^0$ on $\R_+^3\times G_\A^0$
and
$\md\lam_1\md\lam_2\md\ti t^0$ on $\R_+^2\times \ti T_\A^0$,
where the volume of $\ti T_\A^0/\ti T_k$ under $\md\ti t^0$ is $\gC_k^2$.
Moreover, $|\chi(1,\lamb)|=\lam^4$, and so
if $Z^\ast(\Phi,s)$ denotes the zeta function studied in \cite{zfps},
then we have $Z(\Phi,s)=8\Delta_k^{1/2}Z^\ast(\Phi,4s)$.
In \cite{zfps}, it is shown that
$Z^\ast(\Phi,s)$ has a meromorphic continuation to the region $\Re(s)>6$
with only possible singularity in this region at $\Re(s)=8$
with residue
\begin{equation*}
Z_k(2)\gC_k^{-1}\int_{G_{1\A}^0/G_{1k}} dg_{11}^0dg_{12}^0
\cdot\int_{V_\A}\Phi(x)dx.
\end{equation*}
where the measure $dx$ on $V_\A$ is $\Delta_V^{-1/2}$
times that of in this paper.
Thus we arrive at:

\begin{thm}\label{globalanalyticproperties}
Assume that the Schwartz-Bruhat function
$\Phi\in\cS(V_\A)$ has a product form
$\Phi=\otimes_{v\in\gM}\Phi_v$ and
each $\Phi_v\in\cS(V_{k_v})$ is $\cK_v$-invariant.
The zeta function $Z(\Phi,s)$ has a meromorphic continuation to the
region $\re(s)>3/2$ only with a possible simple pole at $s=2$ with residue
\begin{equation*}
\cR_1\prod_{v\in\gM}\int_{V_{k_v}}\Phi_v(x_v)dx_v,
\end{equation*}
where we put
\begin{equation*}
\cR_1=2\Delta_k^{-5/2}\gC_kZ_k(2)^3\prod_{v\in\gMf\cap\gM_\esB}(1-q_v^{-1})^2.
\end{equation*}
\end{thm}
This completes our review of the analytic properties of the global
zeta function.
To arrive at the density theorem from this,
we need various preparations from local theory.
We do it in the next five sections.

\section{The canonical measure on the stabilizer}\label{sec:can}
In this section we shall define a measure on $G_{x\,k_v}^{\circ}$
for $x\in V_{k_v}^{\sst}$ which is canonical in the sense made precise
by Proposition \ref{canonicity}.
Recall that there exists a unique division quaternion algebra $\esB$
up to isomorphism over a local field $F$ other than $\C$, and
that for any separable quadratic extension $L/F$,
there exists a injective homomorphism $L\rightarrow \esB$ of $F$-algebras.
Hence by Proposition \ref{rod},
the set of rational orbits $G_k\backslash V_{k_v}^\sst$ corresponds
to the set of all separable quadratic algebras of $k_v$ if $v\notin\gM_\esB$ and
to the set of all separable quadratic extensions of $k_v$ if $v\in\gM_\esB$.

Following \cite{kayuI}, we
attach to each orbit in $V_{k_v}^{\sst}$ 
where $v\in\gM$, an index or type which
records the arithmetic properties of $v$ and the extension of
$k_v$ corresponding to the orbit. The orbit corresponding to
$k_v\times k_v$ will have the index (ur sp).
(This case does not occur when $v\in\gM_\esB$.)
The orbit corresponding to the unique unramified quadratic extension of
$k_v$ will have the index (rm ur) or (ur ur)
according as $v$ is in $\gM_\esB$ or not.
An orbit corresponding to a ramified quadratic extension of $k_v$
will have the index (rm rm) if
$v\in\gM_{\esB}$ and (ur rm) if $v\notin\gM_{\esB}$.

We first give a normalization 
of the measure on the stabilizer $G_{x\,k_v}^\circ$
for elements of $V_{k_v}^\sst$ of the form $w_u=(1,u)$.
We recall that $k_v[u]$ is isomorphic to
either $k_v\times k_v$ or a quadratic extension of $k_v$
as a $k_v$-algebra.
By using this isomorphism,
we can construct an isomorphism of multiplicative group
\begin{equation*}\label{multiplicativegroup}
k_v[u]^\times\cong\begin{cases}
k_v^{\times}\times k_v^\times &\text{\qquad $w_u$ has type (ur sp),} \\
L_{v,w_u}^\times &\text{\qquad otherwise,} \\
\end{cases}
\end{equation*}
where $L_{v,w_u}$ is the splitting field of $F_{w_u}(v)$ if
this quadratic form is irreducible.
Using the normalized measure of $k_v^\times$ and $L_{v,w_u}^\times$
in Section \ref{sec:not}
(we consider the product measure on $k_v^\times\times k_v^\times$),
we induce a measure $d_u^\times s$ on $k_v[u]^\times$ 
as the pullback measure via the above isomorphism.
We note that this normalization does not depend on
the choice of the isomorphism.

For an element of the form $w_u=(1,u)$, we constructed an isomorphism
\begin{equation*}
\psi_u\colon 
k_v[u]^\times\times k_v[u]^\times
\longrightarrow
G_{w_u\, k_v}^\circ,
\qquad
(s_1,s_2)
\longmapsto
g_{w_u,v}''=(s_1,s_2,g_{s_1s_2})
\end{equation*}
in Section \ref{sec:rev}.
Using this isomorphism and
the product measure $d_u^\times s_1d_u^\times s_2$
on $k_v[u]^\times\times k_v[u]^\times$,
we define a Haar measure $dg_{w_u,v}''$ on $G_{w_u\, k_v}^\circ$ by
\begin{equation*}
dg_{w_u,v}''=(\varphi_u)_\ast(d_u^\times s_1d_u^\times s_2),
\end{equation*}
the pushout measure.
For a general element $x\in V_{k_v}^\sst$
we choose an element $g\in G_{k_v}$ so that
$x=gw_u$ for some $w_u\in V_{k_v}^\sst$,
which is possible by Lemma \ref{lem:orbit}.
Then
\begin{equation*}
i_g\colon
G_{w_u\, k_v}^\circ\longrightarrow G_{x\, k_v}^\circ,
\qquad
g_{w_u,v}''\longmapsto g_{x,v}''=gg_{w_u,v}''g^{-1}
\end{equation*}
gives an isomorphism of groups.
We define the measure $dg_{x,v}''$ on $G_{x\, k_v}^\circ$ by
\begin{equation*}
dg_{x,v}''=(i_g)_\ast(dg_{w_u,v}'').
\end{equation*}
We let $d\ti g_{x,v}''$ on $G_{x\, k_v}^\circ/\ti T_{k_v}$
such that $dg_{x,v}''=d\ti g_{x,v}''\md\ti t_v$.
Note that we defined the measure
$\md \ti t_v$ on $\ti T_{k_v}$ in Section \ref{sec:rev}.

We have to check that these normalizations are well-defined.

\begin{prop}\label{canonicity}
\begin{enumerate}[{\rm (1)}]
\item
The above definition of $dg_{x,v}''$ does not depend on
the choice of $u$ and $g$.
\item Moreover,
suppose that $x,y\in V_{k_v}^{\text{\upshape ss}}$
and that $y=g_{xy}x$ for some
$g_{xy}\in G_{k_v}$. Let $i_{g_{xy}}:G_{y\,k_v}^{\circ}\to
G_{x\,k_v}^{\circ}$ be the isomorphism
$i_{g_{xy}}(g)=g_{xy}^{-1}gg_{xy}$. Then
\begin{equation*}
dg_{y,v}''=i_{g_{xy}}^*(dg_{x,v}'')
\qquad \text{and}\qquad
d\ti g_{y,v}''=i_{g_{xy}}^*(d\ti g_{x,v}'')\,.
\end{equation*}
\end{enumerate}
\end{prop}
\begin{proof}
By the construction of the measures,
a formal consideration shows that it is enough to prove (2)
for $x=(1,u_1),y=(1,u_2)$ where $u_1,u_2\in \esB_{k_v}$.
We write
\begin{equation*}
g_{xy}=(\alpha,\beta,g_2),\quad g_2=\twtw pqrs.
\end{equation*}
Then since $y=g_{xy}x$, we have
\begin{equation}\label{relation}
\beta=(p+qu_1)^{-1}\cdot\alpha^{-1},\quad
u_2=\alpha\cdot \frac{r+su_1}{p+qu_1}\cdot\alpha^{-1}.
\end{equation}
Therefore if we let
\begin{equation*}
\eta\colon
\esB_v\longrightarrow \esB_v,\qquad
\eta(\theta)=\alpha^{-1}\cdot\theta\cdot\alpha,
\end{equation*}
we have $\eta(u_2)=(r+su_1)/(p+qu_1)^{-1}\in k_v[u_1]$
and hence $\eta$ induces
an isomorphism of $k_v$-algebras
\begin{equation}\label{algebra}
\eta\colon k_v[u_2]\longrightarrow k_v[u_1]
\end{equation}
and an isomorphism of groups
\begin{equation}\label{multiplicative}
\eta\colon k_v[u_2]^\times\longrightarrow k_v[u_1]^\times.
\end{equation}
%
Since \eqref{algebra} is an isomorphism of $k_v$-algebras,
\eqref{multiplicative} is a measure preserving map.

Now we show that the diagram
\begin{equation}\label{cd}
\begin{CD}
\left(k_v[u_2]^\times\right)^2
	@>{\psi_{u_2}}>>	G_{y\, k_v}^\circ	\\
@V{(\eta,\eta)}VV		@VV{i_{g_{xy}}}V			\\
\left(k_v[u_1]^\times\right)^2
	@>{\psi_{u_1}}>>	G_{x\, k_v}^\circ	\\
\end{CD}
\end{equation}
is commutative.
Let $s_1,s_2\in k_v[u_2]^\times$.
We compare 
\begin{equation}\label{twoelements}
\psi_{u_1}\circ(\eta,\eta)(s_1,s_2)
\qquad
\text{and}
\qquad
i_{g_{xy}}\circ\psi_{u_2}(s_1,s_2).
\end{equation}
Note that by Proposition \ref{stbstrk},
the $G_{2}$-part of an element of $G_{x}^\circ$
is uniquely determined by its $G_{1}$-part and hence
to prove the above elements are same,
it is enough to verify that their $G_{1}$-parts coincide.
By the definition of the maps, we immediately see
\begin{align*}
\psi_{u_1}\circ(\eta,\eta)(s_1,s_2)
&=(\alpha^{-1}s_1\alpha,\alpha^{-1}s_2\alpha,*),\\
i_{g_{xy}}\circ\psi_{u_2}(s_1,s_2)
&=(\alpha,\beta,g_2)^{-1}(s_1,s_2,*)(\alpha,\beta,g_2)
=(\alpha^{-1}s_1\alpha,\beta s_2\beta^{-1},*).
\end{align*}
Note that we defined $G_{12}$ to be the multiplicative group of the
opposite algebra of $\esB$.
We consider the $G_{12}$-part of the latter element.
By \eqref{relation}, we have
\begin{equation*}
\alpha\beta
=\alpha(p+qu_1)^{-1}\alpha^{-1}
=\eta^{-1}\left((p+qu_1)^{-1}\right)
\in k_v[u_2]
\end{equation*}
and hence commutative with $s_2\in k_v[u_2]$.
Therefore $\alpha\beta s_2=s_2\alpha\beta$ and hence
$\beta s_2\beta^{-1}=\alpha^{-1}s_2\alpha$.
%
This shows that the $G_1$-parts of \eqref{twoelements} coincide
and hence the diagram \eqref{cd} is commutative.
Since 
$(\eta,\eta)\colon 
(k_v[u_2]^\times)^2\rightarrow (k_v[u_1]^\times)^2$
is measure preserving,
the commutativity of the above diagram
establishes the first claim of (2) and the second claim follows from the
observation that $i_{g_{xy}}|_{\ti T_{k_v}}$ is the identity map.
\end{proof}

\section{The local zeta function and the local density}\label{sec:loc}
In this section, 
we make a canonical choice of a measure on the
stabilizer quotient $G_{k_v}/G_{x\,k_v}^{\circ}$
and define the local zeta function.
We also choose a standard orbital representative
for each $G_{k_v}$-orbit in $V_{k_v}$,
and define the the local density $E_v$ for $v\in\gM$
which will show up later in the Euler factor in the density theorem.

We choose a left invariant measure $dg _{x,v}'$ on 
$G_{k_v}/G^{\circ}_{x\,k_v}$ such that 
$dg_v=dg_{x,v}'dg_{x,v}''$.
Recall that we defined invariant measures $dg_v$ and $dg_{x,v}''$
on $G_{k_v}$ and $G_{x\,k_v}^\circ$ in Sections \ref{sec:inv} and
\ref{sec:can}, respectively.
If $g_{xy}\in G_{k_v}$ satisfies
$y=g_{xy}x$ and $i_{g_{xy}}$ is the inner automorphism $g\mapsto
g_{xy}^{-1}gg_{xy}$ of $G_{k_v}$ then
$i_{g_{xy}}(G_{y\,k_v}^{\circ})=G_{x\,k_v}^{\circ}$ and so $i_{g_{xy}}$
induces a homeomorphism $G_{k_v}/ G_{y\,k_v}^{\circ}\to
G_{k_v}/G_{x\,k_v}^{\circ}$,
which we also express by $i_{g_{xy}}$.
\begin{prop}\label{induce}
We have $i_{g_{xy}}^*(dg_{x,v}')=dg_{y,v}'$.
\end{prop}
\begin{proof}
Since the group $G_{k_v}$ is unimodular, 
$i_{g_{xy}}^*(dg_v)=dg_v$.
On the other hand,
we have $i_{g_{xy}}^*(dg_{x,v}'')=dg_{y,v}''$
by Proposition \ref{canonicity}. Hence,
\begin{equation*}
\begin{aligned}
dg_{y,v}' dg''_{y,v}
& = dg_v =i_{g_{xy}}^*(dg_v)\\
& =i_{g_{xy}}^*(dg_{x,v}' dg''_{x,v})
  =i_{g_{xy}}^*(dg_{x,v}')i_{g_{xy}}^*(dg''_{x,v})\\
& =i_{g_{xy}}^*(dg_{x,v}')dg''_{y,v}.
\end{aligned}
\end{equation*}
Therefore $i_{g_{xy}}^*(dg_{x,v}')=dg_{y,v}'$.
\end{proof}

\begin{defn}\label{bxdefinition}
For $v\in\gM$ and $x\in V_{k_v}^{\sst}$ we let $b_{x,v}>0$ be the
constant satisfying the following equation
\begin{equation*}\label{lint}
\int_{G_{k_v}/G^{\circ}_{x\,k_v}} 
f(g_{x,v}'x)\, dg_{x,v}' 
= b_{x,v}\int_{G_{k_v}x}f(y)|P(y)|_v^{-2}\,dy
\end{equation*} 
for any function $f$ on $G_{k_v}x\subset V_{k_v}$
integrable with respect to $dy/|P(y)|_v^2$.
\end{defn}
This is possible because $dy/|P(y)|_v^{2}$ is a
$G_{k_v}$-invariant measure on $V_{k_v}^{\sst}$ and each of the
orbits $G_{k_v}x$ is an open set in $V_{k_v}^{\sst}$.

\begin{prop}\label{bxindep} 
If {\rm $x,y\in V^{\text{ss}}_{k_v}$} and 
$G_{k_v}x=G_{k_v}y$ then $b_{x,v}=b_{y,v}$.  
\end{prop} 
\begin{proof}
Let $f(y)$ be as in Definition \ref{bxdefinition}
and $y=g_{xy}x$ for $g_{xy}\in G_{k_v}$.
Then
\begin{equation*}
\begin{aligned}
\int_{G_{k_v}x}f(y)|P(y)|_v^{-2}\,dy
&= b_{x,v}^{-1}
\int_{G_{k_v}/G^{\circ}_{x\,k_v}} 
f(g_{x,v}'x)\, dg_{x,v}'\\
&= b_{x,v}^{-1}
\int_{G_{k_v}/G^{\circ}_{y\,k_v}} 
f(g_{xy}^{-1}g_{y,v}'y)\, dg_{y,v}'\quad \text{by Proposition \ref{induce}}\\
&= b_{x,v}^{-1}b_{y,v}\int_{G_{k_v}x}f(y)|P(y)|_v^{-2}\,dy.
\end{aligned}
\end{equation*}
Note that the last step is justified
because $dg_{y,v}'$ is left $G_{k_v}$-invariant.
Therefore $b_{x,v}=b_{y,v}$.
\end{proof}

\begin{defn}\label{lzeta}  
For $\Phi\in \cS(V_{k_v})$ and $s\in\C$
we define
\begin{equation*}
\begin{aligned}
Z_{x,v}(\Phi_v,s) 
& =\int_{G_{k_v}/G^{\circ}_{x\,k_v}} 
|\chi(g_{x,v}')|_v^s \Phi_v(g_{x,v}'x)\, dg_{x,v}' \\
\end{aligned}
\end{equation*}
and call it {\em the local zeta function}.
\end{defn}
By the definition of $b_{x,v}$ and 
the equation $P(g_{x,v}'x)=\chi(g_{x,v}')P(x)$, we have
\begin{equation*} \label{oint}
\begin{aligned}
Z_{x,v}(\Phi,s) 
& = \frac{b_{x,v}}{|P(x)|_v^s} \int_{G_{k_v}x} |P(y)|_v^{s-2} \Phi(y)\,dy\,.
\end{aligned}
\end{equation*}
This integral converges absolutely at least when $\re(s)>2$. 
For $x,y\in V_{k_v}^\sst$ lying in the same orbit,
by the above equation and Proposition \ref{bxindep}, we obtain the following.
\begin{prop}\label{sameorbitzeta}
If {\rm $x,y\in V^{\text{ss}}_{k_v}$} and 
$G_{k_v}x=G_{k_v}y$ then
\begin{equation*}
Z_{x,v}(\Phi_v,s)=\frac{|P(y)|_v^s}{|P(x)|_v^s}Z_{y,v}(\Phi_v,s).
\end{equation*}
\end{prop}

By this proposition, we see that the local zeta functions
for the same $G_{k_v}$-orbit are related by a simple equation.
In section \ref{sec:mea}, we define and consider certain Dirichlet series
arising from the global zeta function.
Here, collecting the orbital zeta functions
lying in the same $G_{k_v}$-orbit will be fundamental.
For this purpose, we fix a representative element
for each $G_{k_v}$-orbit in $V_{k_v}^\sst$,
which also has some good arithmetic properties if $v\in\gMf$.

\begin{defn}
For each of $G_{k_v}$-orbits in $V_{k_v}^\sst$,
we choose and fix an element $x$ which satisfies the following condition.
\begin{enumerate}[(1)]
\item If $v\in\gMf$, then $x$ is of the form $(1,u)$
and $u$ generates $\co_{\ti k_v(x)}$ over $\co_v$
via the identification $k_v[u]\cong \ti k_v(x)$.
\item If $v\in\gMi$, then $|P(x)|_v=1$.
\end{enumerate}
We call such fixed orbital representatives as
the {\em standard orbital representatives}.
\end{defn}

If $v\in\gMf$,
for any standard representative $x=(1,u)\in V_{k_v}^\sst$,
$u$ is a root of $F_x(v_1,-1)$
and so the discriminant $P(x)$ of $F_x(v)$
generates the ideal $\Delta_{k_v(x)/k_v}$.

\begin{defn}
For any $v\in\gMf$, let $\Phi_{v,0}$ be
the characteristic function of $V_{\co_v}$.
Also we put
\begin{equation*}
Z_{x,v}(s)=Z_{x,v}(\Phi_{v,0},s).
\end{equation*}
We call $Z_{x,v}(s)$ for any standard orbital representative $x$
a {\em standard local zeta function} of $x$.
\end{defn}

To describe estimates of Dirichlet series,
we introduce the following notation.
\begin{defn}
Suppose that we have Dirichlet series
$L_i(s)=\sum_{m=1}^{\infty}\ell_{i,m}m^{-s}$ for $i=1,2$. If
$\ell_{1,m}\leq\ell_{2,m}$ for all $m\geq1$ then we shall write
$L_1(s)\cleq L_2(s)$.
\end{defn}

We set $S_0=\gMi\cup\gMdy\cup\gM_\esB$.
To carry out the filtering process,
we need a uniform estimate of the standard local zeta functions.
The following proposition concerning the standard local zeta functions
for $v\notin S_0$ is proved in
\cite[Corollary 8.24, Proposition 9.25]{kayuI}.
Since $S_0$ is a finite set, the result is enough for our purposes.

\begin{prop}\label{estimate}
Let $v\notin S_0$ and $x\in V_{k_v}^\sst$
be one of the standard representatives.
Then
$Z_{x,v}(s)$ can be expressed as
\begin{equation*}
Z_{x,v}(s)=\sum_{n\geq0}\frac{a_{x,v,n}}{q_v^{ns}}
\end{equation*}
with $a_{x,v,0}=1$ and $a_{x,v,n}\geq0$ for all $n$.
Also let us define
\begin{equation*}
L_v(s) =
\frac{1+29q_v^{-2(s-1)}-21q_v^{-4(s-1)}+7q_v^{-6(s-1)}} 
{(1-q_v^{-(2s-1)})(1-q_v^{-2(s-1)})^4}.
\end{equation*}
Then $Z_{x,v}(s)\cleq L_v(s)$.
\end{prop}

Now we define the local density.

\begin{defn}
Assume $x\in V_{k_v}^\sst$ is a standard orbital representative.
We define
\begin{equation*}
\varepsilon_v(x)=\frac{|P(x)|_v^2}{b_{x,v}}.
\end{equation*}
Also we define the {\em local density} at $v$ by
\begin{equation*}
E_v=\sum_x \varepsilon_v(x)
\end{equation*}
where the sum is over all standard representatives
for orbits in $G_{k_v}\backslash V_{k_v}^\sst$.
\end{defn}

These values plays an essential role in the density theorem.
The purpose in the next three sections are
to compute the local densities.
To make the density theorem more precise,
it is better to evaluate $\varepsilon_v(x)$ separately
rather than the sum $E_v$.
We compute for $v\in\gMf$ in Sections \ref{sec:fin}, \ref{sec:ram}
and for $v\in\gMi$ in Section \ref{sec:inf}.
For $v\notin\gM_\esB$, those were already almost carried out in
\cite{kayuI,kayuII} and 
except for a refinement
for dyadic places in Proposition \ref{efur},
we quote their result.

\begin{rem}\label{rem:kayu}
We briefly compare the definition of standard orbital representatives
and the value $\varepsilon_v(x)$
in \cite{kayuI} and in this paper for $v\notin\gM_\esB$,
to confirm that we can directly use their result.
Let $G'_{k_v}$ denote the group
of the representation $V_{k_v}$ used in \cite{kayuI}.
Then one can easily see that the isomorphism $G'_{k_v}\rightarrow G_{k_v}$
given by $(g_1,g_2,g_3)\mapsto (g_1,{}^tg_2,g_3)$ is compatible
with their actions on $V_{k_v}$.
If we identify these groups using the isomorphism,
we immediately see that our choice of measure on $G_{k_v}$
coincides to that of in \cite{kayuI},
and moreover, measures on $G_{x\, k_v}^\circ$ also.
The latter claim holds because both papers used the isomorphism in
Proposition \ref{stbstrk} to normalize the measures on $G_{x\, k_v}^\circ$.
The normalization of the measures on $G_{k_v}/G_{x\,k_v}^\circ$
are slightly different, but from the definitions
we can easily see that
the constants $b_{x,v}$'s coincide.
Although our choice of the standard orbital representatives $x$
for $v\in\gMf$ is also slightly different,
the values of $|P(x)|_v$ coincide since
the standard orbital representative $x$ in \cite{kayuI}
are also chosen so that $P(x)$ generate $\Delta_{k_v(x)/k_v}$.
Since $\varepsilon_v(x)$ is determined only by $|P(x)|_v$ and $b_{x,v}$,
this observation shows that our $\varepsilon_v(x)$'s coincide to
those of \cite{kayuI}.
\end{rem}

\section{Computation of the local densities at finite unramified places}
\label{sec:fin}

In this and next sections, we assume $v\in\gMf$.
We first introduce some notations for these sections.
For any $v\in\gMf$ we shall put $2\co_v=\gp_v^{m_v}$.
Of course $m_v=0$ unless $v\in\gMdy$.
If $x\in V_{k_v}^{\sst}$ then let
$\Delta_{k_v(x)/k_v}=\gp_v^{\delta_{x,v}}$. 
It is well-known that if $k_v(x)/k_v$ is ramified
then $\del_{x,v}$ takes one of the values $2,4,\dots,2m_v,2m_v+1$.
(In the case $v\notin\gMdy$ and hence $m_v=0$,
this should be counted as
$\delta_x$ only takes the value $1$.)

We now assume $v\notin\gM_\esB$.
The following propositions are proved in
\cite[Lemma 7.3]{kayuI} and
\cite[Propositions 4.14, 4.15, 4.25]{kayuII}.

\begin{prop}\label{efuu}
Assume $v\notin\gM_\esB$.
Let $x\in V_{k_v}^\sst$ be one of the standard representative.
\begin{enumerate}[{\rm (1)}]
\item If $x$ has type {\upshape (ur~sp)} then 
$\varepsilon_v(x)=2^{-1}(1+q_v^{-1})(1-q_v^{-2})^2$.
\item If $x$ has type {\upshape (ur~ur)} then
$\varepsilon_v(x)=2^{-1}(1-q_v^{-1})^3(1-q_v^{-2})$.
\end{enumerate}
\end{prop}

\begin{prop}\label{efurKY}
Assume $v\notin\gM_\esB$.
Let $x\in V_{k_v}^\sst$ be one of the standard representative.
\begin{enumerate}[{\rm (1)}]
\item If $v\notin\gMdy$ and 
$x$ has type {\upshape (ur~rm)} then
$\varepsilon_v(x)=2^{-1}q_v^{-1}(1-q_v^{-1})(1-q_v^{-2})^3$.
\item
If $v\in\gMdy$ then
\begin{equation*}
\begin{aligned}
\sum_{2\leq \del_{x,v}= 2\ell\leq 2m_v}\varepsilon_v(x) & = 
(1-q_v^{-1})^2(1-q_v^{-2})^3q_v^{-\ell}, \\
\sum_{\del_{x,v}= 2m_v+1}\varepsilon_v(x) & = 
(1-q_v^{-1})(1-q_v^{-2})^3q_v^{-(m_v+1)},
\end{aligned}
\end{equation*}
where $x$ runs through all the standard representative
with the given condition of discriminants.
\end{enumerate}
\end{prop}

\begin{prop}\label{ldsf}
Let $v\notin\gM_\esB$. Then
\begin{equation*}
E_v=(1-q_v^{-2})(1-3q_v^{-3}+2q_v^{-4}+q_v^{-5}-q_v^{-6}).
\end{equation*}
\end{prop}

These results are already enough to prove our density theorems.
However, if we could know the value $\varepsilon_v(x)$
for $v\in\gMdy$ in the Proposition \ref{efurKY} solely,
then the density theorems become finer.
In this section we refine Proposition \ref{efurKY}
to the following.

\begin{prop}\label{efur}
Assume $v\notin\gM_\esB$.
Let $x\in V_{k_v}^\sst$ be a standard representative
with the type {\upshape (ur~rm)}.
Then
\begin{equation*}
\varepsilon_v(x)
=	2^{-1}|\Delta_{k_v(x)/k_v}|_v^{-1}(1-q_v^{-1})(1-q_v^{-2})^3.
\end{equation*}
\end{prop}

It is well known that there are
$2q_v^{l-1}(q_v-1)$
numbers of quadratic extensions of $k_v$
with the absolute value of the relative discriminant $q_v^{2l}$
if $1\leq l\leq m_v$ and
$2q_v^{m_v}$
numbers of quadratic extensions of $k_v$
with the absolute value of the relative discriminant $q_v^{2m_v+1}$.
Hence this is in fact a refinement of Proposition \ref{efurKY}.
We give the proof of this proposition after we prove Lemma \ref{stabur}.

Let $L/k_v$ be a quadratic ramified extension,
$\varpi$ a uniformizer of $L$,
and $\varpi^\tau$ the conjugate of $\varpi$ with respect to $L/k_v$,
henceforth fixed.
We put $a_1=\varpi+\varpi^\tau, a_2=\varpi\varpi^\tau$.
Following \cite{kayuI}, we let
\begin{equation}\label{ystn}
x=(x_1,x_2),\quad\text{where}\quad
x_1=\twtw 011{a_1}\ \text{and}\ 
x_2=\twtw 1{a_1}{a_1}{a_1^2-a_2}.
\end{equation}
Then
$F_x(v)
=	-(v_1^2+a_1v_1v_2+a_2v_2^2)
	=-(v_1+\varpi v_2)(v_1+\varpi^\tau v_2)$
and hence $L\cong k_v(x)$ and
$P(x)$ generates the ideal $\Delta_{k_v(x)/k_v}$.
Therefore
we can replace the standard representative
for the orbit corresponding to $L$
to this $x$ to compute $\varepsilon_v(x)=|P(x)|_v^2b_{x,v}^{-1}$.

The following lemma is a consequence of
\cite[Lemma 7.3]{kayuI} and \cite[Proposition 3.2]{kayuII}.

\begin{lem}\label{evur}
We have $\varepsilon_v(x)=\vol(\cK_vx)$.
\end{lem}
We compute $\vol(\cK_vx)$ with a slight
modification of the method in \cite{kayuII},
along the line of \cite{pcbq}.
To begin with we introduce some notations,
which we also use to consider similar problems in Section \ref{sec:ram}.
We regard $\cK_v$ as the set of $\co_v$-rational points $G_{\co_v}$
of a group scheme $G=\gl(2)\times\gl(2)\times\gl(2)$ defined over $\co_v$
acting on a module scheme
$V={\rm M}(2,2)\oplus {\rm M}(2,2)$ also defined over $\co_v$.
Then since $x$ is an $\co_v$-rational point of $V$,
we can consider the stabilizer of $x$
as a group scheme also defined over $\co_v$,
in the sense of \cite{mufo}.
Let $G_x$ denote this group scheme.
Note that this definition of $G_x$ differs from \cite{kayuII, kayuIII}.
Let $i$ be a positive integer.
For an $\co_v$-scheme $X$, let $\prj_{X,i}$ denote the reduction map
$X_{\co_v}\rightarrow X_{\co_v/\gp_v^i}$.
If the situation is obvious we drop $X$ and write $\prj_i$ instead.
For rational points $y_1,y_2\in X_{\co_v}$,
we use the notation $y_1\equiv y_2\;(\gp_v^i)$ if
$\prj_i(y_1)=\prj_i(y_2)$.
We also use the notation ``$y$ mod $\gp_v^i$'' for $\prj_i(y)$.

For the element $x$ of the form \eqref{ystn}, let
\begin{equation*}
A_x(c,d)=\twtw cd{-a_2d}{c+a_1d}.
\end{equation*}
Then if $A_x(c_i,d_i)\in\gl(2)_{\co_v}$ for $i=1,2$,
by computation we could see that the element
\begin{equation*}
(A_x(c_1,d_1),A_x(c_2,d_2),A_x(c_1,d_1)^{-1}A_x(c_2,d_2)^{-1})\in \cK_v
\end{equation*}
stabilizes $x$.
Let $N_{x\,\co_v}$ denote the subgroup of $\cK_v$
consisting of elements of the form above.
We naturally regard $N_{x\,\co_v}$
as the set of $\co_v$-rational points
of a group scheme $N_x$,
which is a subgroup of $G_x$,
defined over $\co_v$.

\begin{prop}\label{stbstru}
We have $N_x\cong (\co_{k_v(x)}^\times)^2$ as a group scheme over $\co_v$.
\end{prop}
\begin{proof}
Let $R$ be any $\co_v$-algebra.
Then we could see that the map
\begin{equation*}
(A_x(c_1,d_1),A_x(c_2,d_2),A_x(c_1,d_1)^{-1}A_x(c_2,d_2)^{-1})
\mapsto (c_1+\varpi d_1,c_2+\varpi d_2)
\end{equation*}
gives an isomorphism between $N_{x\, R}$ and
$\{(\co_{k_v(x)}\otimes_{\co_v} R)^\times\}^2$,
and this map, denoted by $\psi_{x,R}$,
satisfies the usual functorial property
with respect to homomorphism of $\co_v$-algebras.
This shows that there exists an isomorphism
$\psi_x\colon N_x\rightarrow (\co_{k_v(x)}^\times)^2$
as groups schemes over $\co_v$
such that $\psi_{x,R}$ is the induced isomorphism for all $R$.
\end{proof}

We now consider the orbit $\cK_vx$.
The approach in \cite{kayuII} is
to consider modulo $\gp_v$ congruence condition on $V_{\co_v}$
to compute the sum $\sum_{x}\vol(\cK_vx)$ where $x$ runs through
all the standard representatives with the given relative discriminant.
Let $n=\delta_{x,v}+2m_v+1$ as in \cite{pcbq}.
Then, as we demonstrate below,
deliberation of the congruence relation of modulo $\gp_v^n$
allows us to treat the orbit $\cK_vx$ solely.
We note that the idea of considering modulo a certain high power of prime ideal
is already presented in \cite{kayuII} and used to compute $\varepsilon_v(x)$
in some other cases.

\begin{defn}
We define $\cD=\{y\in V_{\co_v}\mid y\equiv x\,(\gp_v^n)\}$.
\end{defn}
\begin{lem}\label{Dcontainedur}
We have $\cD\subset \cK_v x$.
\end{lem}
\begin{proof}
Let $y\in \cD$. First we show $y\in G_{k_v}x$.
Since $P(y)\equiv P(x)\;(\gp_v^n)$
and $\ord_v(P(x))=\delta_{x,v}$,
we have $P(y)/P(x)\equiv 1\;(\gp_v^{2m_v+1})$.
Then by Hensel's lemma,
we have $P(y)/P(x)\in (k_v^\times)^2$.
Therefore the splitting fields of $F_x(v)$ and $F_y(v)$ coincide
and hence by Lemma \ref{rod}, we have $y\in G_{k_v}x$.
The rest of argument is exactly the same as
that of \cite{kayuI, kayuII} and we omit it.
\end{proof}

\begin{lem}\label{stabur}
We have $[G_{x\,\co_v/\gp_v^n}:N_{x\,\co_v/\gp_v^n}]=2q_v^{\delta_{x,v}}$.
\end{lem}
\begin{proof}
The same argument as in the proof of \cite[Proposition 4.15]{kayuII}
shows that each right coset space of 
$N_{x\,\co_v/\gp_v^n}\backslash G_{x\,\co_v/\gp_v^n}$
contains exactly one element of the form
$g=(g_1,g_2)$, $g_1=(g_{11},g_{12})$ with
\begin{equation*}
g_{11}=\twtw 10us, g_{22}=\twtw 1v0t, g_2=\twtw\alpha\beta\gamma\delta.
\end{equation*}
Hence we will consider when
such an element actually lies in $G_{x\,\co_v/\gp_v^n}$.
Suppose that $g$ is in the form above and
$gx=x$ in $V_{\co_v/\gp_v^n}$.
We put $y=(y_1,y_2)=(g_1,1)x$.
Then by computation we have
\begin{equation*}
y_1=\twtw 0ts\ast,\quad
y_2=\twtw 1{v+a_1t}{u+a_1s}\ast.
\end{equation*}
Therefore, by comparing the (1,1), (1,2) and (2,1)-entries of
$x_1$ and $\alpha y_1+\beta y_2$,
we have $\beta=0$ and $s=t=\alpha^{-1}$.
Also under the condition $s=t$, from $x_1=\gamma y_1+\delta y_2$
we have $\delta=1$, $u=v$ and $\gamma=s^{-1}(a_1-u-a_1s)$.
Under these equations, we have
\begin{align*}
\alpha y_1+\beta y_2	&=\twtw 011{2u+a_1s}\\
\gamma y_1+\delta y_2	&=\twtw 1{a_1}{a_1}{-u^2+(2u+a_1s)a_1-a_1su-a_2s^2}
\end{align*}
and therefore we could see that $gx=x$ if and only if
\begin{equation*}
2u+a_1s=a_1\quad\text{and}\quad
u^2+a_1su+a_2s^2=a_2\quad \text{in}\ \co_v/\gp_v^n.
\end{equation*}
This system is exactly the same as that we considered
in \cite[Lemma 4.7]{pcbq} and it has $2q_v^{\delta_{x,v}}$
solutions in all.
\end{proof}

We are now ready to prove Proposition \ref{efur}.
Let $\gr_n$ be the reduction map $G_{\co_v}\rightarrow G_{\co_v/\gp_v^n}$.
Then by Lemma \ref{Dcontainedur}, the set $\cK_v x=G_{\co_v}x$ is equal to
$\#(G_{\co_v}/\gr_n^{-1}(G_{x\,\co_v/\gp_v^n}))$ number of
disjoint copies of $\cD$.
Since
\begin{equation*}
G_{\co_v}/\gr_n^{-1}(G_{x\,\co_v/\gp_v^n})\cong
G_{\co_v/\gp_v^n}/G_{x\,\co_v/\gp_v^n},
\end{equation*}
by Lemma \ref{stabur} we have
\begin{equation*}
\begin{aligned}
\vol(\cK_v x)
&=\vol(\cD)\cdot\frac{\#(G_{\co_v/\gp_v^n})}
{2q_v^{\delta_x}\cdot\#(N_{x\,\co_v/\gp_v^n})}
=q_v^{-8n}\cdot
\frac{\left\{q_v^{4n}(1-q_v^{-1})(1-q_v^{-2})\right\}^3}
{2q_v^{\delta_x}\cdot \left\{q_v^{2n}(1-q_v^{-1})\right\}^2}\\
&=2^{-1}q_v^{-\delta_{x}}(1-q_v^{-1})(1-q_v^{-2})^3.
\end{aligned}
\end{equation*}
Since $|\Delta_{k(x)/k}|_v=q_v^{\delta_x}$, we obtained the desired result.

\section{Computation of the local densities at finite ramified places}
\label{sec:ram}
In this section we assume $v\in\gM_\esB$ and so
$\esB_v$ is a non-split quaternion algebra of $k_v$.
We briefly recall the algebraic structure of $\esB_v$
and prepare the notations to begin with.
We take a commutative subalgebra $F_v$ of $\esB_v$ so that
$F_v$ is a quadratic unramified extension of $k_v$
and henceforth fixed in this section.
Let $\sigma$ denote the non-trivial element of ${\rm Gal}(F_v/k_v)$.
Then for any prime element $\pi_v\in k_v$,
$\esB_v$ can be identified with $F_v\oplus F_v\sqrt{\pi_v}$
as a left vector space of $F_v$ and the multiplication law is
given by $\sqrt{\pi_v}\alpha=\alpha^\sigma\sqrt{\pi_v}$
for $\alpha\in F_v$.
For $a\in \esB_v$, let $a^\ast$ be its involution.
Then for $\alpha,\beta\in F_v$,
$(\alpha+\beta\sqrt{\pi_v})^\ast=\alpha^\sigma-\beta\sqrt{\pi_v}$.
Hence the reduced trace $\trc$ and the reduced norm $\nr$ of $\esB_v$ is given by
\begin{equation*}
\trc(\alpha+\beta\sqrt{\pi_v})=\alpha+\alpha^\sigma,
\quad
\nr(\alpha+\beta\sqrt{\pi_v})=\alpha\alpha^\sigma-\pi_v\beta\beta^\sigma,
\end{equation*}
for $\alpha,\beta\in F_v$.
The map $u\mapsto\ord_v(\nr(u))$ defines a discrete valuation of
$\esB_v$, and it is well known that
$\co_{\esB_v}=\{u\mid |\nr(u)|_v\leq1\},\co_{\esB_v}^\times=\{u\mid |\nr(u)|_v=1\}$.
If we restrict the reduced norm
to any quadratic subfield $L_v$,
it coincides with the norm map 
$\n_{L_v/k_v}$ of the extension $L_v/k_v$.
Hence $\co_{\esB_v}\cap L_v=\co_{L_v}$
and $\co_{\esB_v}^\times \cap L_v=\co_{L_v}^\times$. 
We fix an element $\theta\in\co_{F_v}$ 
so that $\co_{F_v}=\co_v[\theta]$.
By computation we have the following.
\begin{lem}\label{integralring}
We have
\begin{equation*}
\co_{\esB_v}=\{\alpha+\beta\sqrt{\pi_v}\mid \alpha,\beta\in\co_{F_v}\},\quad
\co_{\esB_v}^\times=\{\alpha+\beta\sqrt{\pi_v}\mid \alpha\in\co_{F_v}^\times,
\beta\in\co_{F_v}\}.
\end{equation*}
For any quadratic ramified extension $L_v$ of $k_v$
contained in $\esB_v$,
we can also write
\begin{equation*}
\co_{\esB_v}=\{\alpha+\beta\theta\mid \alpha,\beta\in\co_{L_v}\},\quad
\co_{\esB_v}^\times=\{\alpha+\beta\theta\mid \alpha,\beta\in\co_{L_v},
\alpha\in\co_{L_v}^\times\ \text{or}\ \beta\in\co_{L_v}^\times\},
\end{equation*}
and moreover, by changing $\theta$ if necessary,
the multiplication law is given by $\theta\alpha=\alpha^\tau\theta$
for $\alpha\in L_v$ where $\tau$ denote the non-trivial element of
$\gal(L_v/k_v)$.
\end{lem}

As in Section \ref{sec:fin}, we can and shall regard
$\cK_v$
as the set of $\co_v$-rational points $G_{\co_v}$ 
of a group scheme $G$ defined over $\co_v$
acting on a module scheme $V$ also defined over $\co_v$.
For example,
the group
$G_{\co_v}=
\co_{\esB_v}^\times\times(\co_{\esB_v})^\times\times\gl(2)_{\co_v}$
acts on the module $V_{\co_v}=\co_{\esB_v}\oplus\co_{\esB_v}$.
Then any standard orbital representative $x$ is an element of $V_{\co_v}$
and as in Section \ref{sec:fin},
we regard the stabilizer $G_x$
as a group scheme defined over $\co_v$.
If $x\in V_{\co_v}$, then $F_x(v)\in\sym^2\co_v^2$.
We also regard $\sym^2\co_v^2$ as a module scheme over $\co_v$
and the map $x\mapsto F_x(v)$ as a morphism of schemes.
We continue to use the notation $\prj_i$ defined in Section \ref{sec:fin}.
Further, for any quadratic extension $L_v$ of $k_v$,
we use the abbreviation $\co_{L_v}/\gp_v^i$ for $\prj_i(\co_{L_v})$,
where we regard $\co_{L_v}$ as a scheme over $\co_v$.
For example,
\begin{equation*}
\co_{F_v}/\gp_v^i=(\co_v/\gp_v^i)[\theta]=
\{\alpha+\beta\theta\mid \alpha,\beta\in\co_v/\gp_v^i\}.
\end{equation*}
In this section, we will express $g\in G_{k_v}$ as
\begin{equation*}
g=(g_{11},g_{12},g_2),\quad
g_2=\twtw pqrs.
\end{equation*}

\begin{prop}\label{stbstr}
Let $x$ be one of the standard representative.
Then there exist an injective homomorphism
$(\co_{\ti k_v(x)}^\times)^2\rightarrow G_x$
as a group scheme over $\co_v$.
\end{prop}
\begin{proof}
Let $x=(1,u)$.
We construct the
injective homomorphism
\begin{equation*}
\psi_{u\, R}\colon
\{(\co_{\ti k_v(x)}\otimes R)^\times\}^2\longrightarrow G_{x\, R}
\end{equation*}
for any commutative $\co_v$-algebra $R$.

We put $\ti R(x)= \co_{\ti k_v(x)}\otimes R$.
Note that $\ti R(x)=R[u]$ is a subalgebra of
$\co_{\esB_v}\otimes R$ and is commutative.
Since $\{1,u\}$ is a $\co_v$-basis of $\co_v[u]$,
this is also an $R$-basis of $\ti R(x)$.
Let $s_1,s_2\in\ti R(x)^\times$.
Then $\{s_1s_2,s_1s_2u\}$ is also an $R$-basis of $\ti R(x)$,
and so there exists a unique element $g=g_{s_1s_2}\in\gl(2)_R$ such that
$g\,\,^{t}(s_1s_2,s_1s_2u)={}^t(1,u)$.
Hence
\begin{equation*}
\psi_{u\, R}\colon (s_1,s_2)\longmapsto (s_1,s_2,g_{s_1s_2})
\end{equation*}
gives an injective homomorphism
from $(\ti R(x)^\times)^2$ to $G_{xR}$,
and as in the proof of Proposition \ref{stbstru},
we can regard this map as the induced one from the morphism of schemes.
\end{proof}
Let $N_x$ denote the image of this homomorphism,
which is a subgroup of $G_x$.

\begin{prop}\label{maxstab}
Let $x\in V_{k_v}^\sst$
be one of the standard representatives.
Then
\begin{equation*}
\int_{\cK_v\cap G_{x\,k_v}^\circ}dg_{x,v}''=1.
\end{equation*}
\end{prop}
\begin{proof}
Let $x=(1,u)$ be a standard representative.
We claim that $\psi_u^{-1}(\cK_v\cap G_{x\,k_v}^\circ)=(\co_v[u]^\times)^2$
where $\psi_u$ is defined in Section \ref{sec:can}.
The inclusion 
$\psi_u^{-1}(\cK_v\cap G_{x\,k_v}^\circ)\subset(\co_{k_v[u]}^\times)^2$
follows from $\co_{\esB_v}^\times\cap k_v(x)=\co_{k_v(x)}^\times$.
Let $s_1,s_2$ be elements of $\co_{k_v[u]}^\times$.
Then since $\{s_1s_2,s_1s_2u\}$ also forms a $\co_v$-basis of $\co_{k_v[u]}$,
we have $g_{s_1s_2}\in\gl(2)_{\co_v}$.
This shows the reverse inclusion.
Now the proposition follows from the definition of $dg_{x,v}''$.
\end{proof}

The following simple observation will be sometimes useful 
in the concrete calculations below.
This easily follows from Proposition \ref{stbstrk}
and the properties of the norm map of the quadratic extension
of local fields.
\begin{lem}\label{orderstab}
We define
\begin{equation*}
\varsigma:G_{k_v}\longrightarrow \Z^2,
\quad
{\text as}
\quad
g\longmapsto (\ord_v(\nr(g_{11})),\ord_v(\nr(g_{12}))).
\end{equation*}
Then the image $\varsigma(G_{x\, k_v}^\circ)$ is
$(2\Z)^2$ if $x$ corresponds to the quadratic unramified extension
and $\Z^2$ if $x$ corresponds to a quadratic ramified extension.
\end{lem}

From now on we consider the case $k_v(x)$ is unramified
and ramified separately.
We first consider the former case.
Till Proposition \ref{efru},
we assume $x$ has type (rm ur).
We note that in this case the polynomial 
$(F_x(v)\mod\gp_v)\in\sym^2(\co_v/\gp_v)^2$
is irreducible and especially $F_x(0,1)\in\co_v^\times$.
By changing the choice of the included unramified extension
$F_v$ and the generator of the integer ring $\theta$ if necessary,
we may assume $x=(1,\theta)$.
Let us write $\theta=a+b\theta^\sigma,a\in\co_v,b\in\co_v^\times$
and set $\tau_\theta=\twtw 10ab\in\gl(2)_{\co_v}$.
We fix a prime element $\pi_v\in k_v$ and
put $\tau_x=(\sqrt{\pi_v}^{-1},\sqrt{\pi_v},\tau_\theta)$,
which then generates the non-trivial class of $G_{x\, k_v}/G_{x\, k_v}^\circ$.
\begin{lem}\label{taux}
Let $x$ have type {\upshape (rm ur)}.
Then
$\cK_vG_{x\,k_v}=\cK_vG_{x\,k_v}^\circ\amalg\tau_x \cK_vG_{x\,k_v}^\circ$.
\end{lem}
\begin{proof}
Since $\ord_v(\nr(\sqrt{\pi_v}^{-1}))=\ord_v(\pi_v^{-1})=-1$,
we have $\tau_x\notin \cK_vG_{x\,k_v}^\circ$
as a consequence of Lemma \ref{orderstab}.
Now the lemma follows since 
$\tau_x$ is a normalizer of the group $\cK_v$.
\end{proof}
\begin{lem}\label{evru}
Let $x$ have type {\upshape (rm ur)}.
Then $\varepsilon_v(x)=2^{-1}\vol(\cK_vx)$.
\end{lem}
\begin{proof}
By the definition of $dg_{x,v}'$, Proposition \ref{maxstab}
and Lemma \ref{taux},
\begin{align*}
1=\int_{\cK_v}dg_v
&=\int_{\cK_v G_{x\,k_v}^{\circ}/G_{x\,k_v}^{\circ}}dg_{x,v}'
\cdot
\int_{\cK_v\cap G_{x\,k_v}^{\circ}} dg_{x,v}''\\
&=\int_{\cK_v G_{x\,k_v}^{\circ}/G_{x\,k_v}^{\circ}}dg_{x,v}'\\
&=\frac12\int_{\cK_v G_{x\,k_v}/G_{x\,k_v}^{\circ}}dg_{x,v}'.
\end{align*}
Hence, if we let $\Phi_v$ be the characteristic function of $\cK_vx$,
by Definition \ref{bxdefinition} we have
\begin{align*}
2=\int_{\cK_v G_{x\,k_v}/G_{x\,k_v}^{\circ}}dg_{x,v}'
&=\int_{G_{k_v}/G_{x\,k_v}^{\circ}}\Phi_v(g_{x,v}'x)\,dg_{x,v}'
=b_{x,v}\int_{G_{k_v}x}\Phi_v(y)|P(y)|^{-2}_v\,dy\\
&=b_{x,v}\int_{\cK_vx}|P(y)|_v^{-2}\,dy.
\end{align*}
Since $|P(y)|_v=|P(x)|_v$ for all $y\in \cK_vx$, we have
$\varepsilon_v(x)=|P(x)|_v^2b_{x,v}^{-1}=2^{-1}\vol(\cK_vx)$.
\end{proof}
We will compute $\vol(\cK_vx)$.
In the case $k_v(x)$ is unramified extension,
it is enough to consider the congruence relation of modulo $\gp_v$.
\begin{defn}
We define $\cD=\{y\in V_{\co_v}\mid y\equiv x\,(\gp_v)\}$.
\end{defn}
\begin{lem}\label{Dcontainedru}
We have $\cD\subset \cK_vx$.
\end{lem}
\begin{proof}
Let $y\in\cD$. Since
$(F_y(v)\mod \gp_v)= (F_x(v)\mod \gp_v)\in\sym^2(\co_v/\gp_v)^2$,
the splitting field of $F_y(v)$ is the quadratic unramified extension.
Hence, $y\in G_{k_v}x$.
Let $y=gx, g=(g_{11},g_{12},g_2)\in G_{k_v}$.
Note that
\begin{equation*}
|\chi(g)|_v=|\nr(g_{11})\nr(g_{12})\det(g_2)|_v^2=1
\end{equation*}
since $|P(y)|_v=|P(x)|_v$.
We will show that $g\in \cK_vG_{x\,k_v}$.
By Lemma \ref{orderstab},
multiplying an element of $G_{x\,k_v}^\circ$
and $\tau_x$ if necessary, we may assume that
$g$ satisfies either one of the following conditions.
\begin{enumerate}
\item[(A)] $|\nr(g_{11})|_v=|\nr(g_{12})|_v=1$, (hence $|\det(g_2)|_v=1$.)
\item[(B)] $|\nr(g_{11})|_v=q_v,|\nr(g_{12})|_v=1$,
(hence $|\det(g_2)|_v=q_v^{-1}$.)
\end{enumerate}
From the definition of the representation we have
\begin{equation*}
F_y(v)=\nr(g_{11})\nr(g_{12})F_x(vg_2)
\end{equation*}
and hence
\begin{align*}
F_y(1,0)&=\nr(g_{11})\nr(g_{12})N_{F_v/k_v}(p+q\theta),\\
F_y(0,1)&=\nr(g_{11})\nr(g_{12})N_{F_v/k_v}(r+s\theta).
\end{align*}
On the other side, since $F_x(v)\equiv F_y(v)\;(\gp_v)$,
both $F_y(1,0)$ and $F_y(0,1)$ are units of $\co_v$.
If $g$ satisfies the condition (B), then
$\ord_v(\n_{F_v/k_v}(p+q\theta))$ must be $1$.
But this is a contradiction since $F_v/k_v$ is
the quadratic unramified extension.
Hence we assume $g$ satisfies the condition (A).
Then both $\n_{F_v/k_v}(p+q\theta)$ and $\n_{F_v/k_v}(r+s\theta)$
are elements of $\co_v^\times$ and so $p,q,r,s\in\co_v$.
Since $|\det(g_2)|_v=1$, we conclude $g_2\in\gl(2)_{\co_v}$.
Thus $g\in \cK_v$ and the lemma follows.
\end{proof}
\begin{lem}\label{stabru}
We have $G_{x\,\co_v/\gp_v}=N_{x\,\co_v/\gp_v}$.
\end{lem}
\begin{proof}
In the proof of this lemma, if we have $y\equiv y'\;(\gp_v)$
for any two $\co_v$-rational points of an $\co_v$-scheme,
we drop $(\gp_v)$ and simply write $y\equiv y'$ instead.
Clearly $G_{x\,\co_v/\gp_v}\supset N_{x\,\co_v/\gp_v}$
and hence we prove the reverse inclusion.
Let $g=(g_{11},g_{12},g_2)\in G_{x\,\co_v/\gp_v}$.
We choose representatives of $g_{11},g_{12},g_2$
in $G_{11k_v},G_{12k_v},G_{2k_v}$ and use the same notation for them.
By Lemma \ref{integralring} and Proposition \ref{stbstr},
multiplying an element of $N_{x\,\co_v/\gp_v}$ if necessary,
we assume that
\begin{equation*}
g_{11}=1+\alpha\sqrt{\pi_v},\quad
g_{12}=1+\beta\sqrt{\pi_v},
\end{equation*}
where $\alpha,\beta\in\co_{F_v}$.
Put $y=(y_1,y_2)=(g_{11},g_{12},1)x$.
Then by computation we have
\begin{equation*}
y_1\equiv1+(\alpha+\beta)\sqrt{\pi_v},\quad
y_2\equiv\theta+(\alpha\theta^\sigma+\beta\theta)\sqrt{\pi_v}.
\end{equation*}
Since ${}^t(y_1,y_2)\equiv g_2^{-1}\,{}^t(1,\theta)$
and $1,\theta\in\co_{F_v}$,
we have
$\prj_1(y_1),\prj_1(y_2)\in\co_{F_v}/\gp_v$.
Hence $\alpha+\beta\equiv0,\alpha\theta^\sigma+\beta\theta\equiv0$.
Then since $\theta-\theta^\sigma\in\co_{F_v}^\times$,
we have $\alpha\equiv\beta\equiv0$ and hence 
$g_{11}\equiv g_{12}\equiv1$.
Then $g_2\equiv I_2$ and this shows
$G_{x\,\co_v/\gp_v}\subset N_{x\,\co_v/\gp_v}$.
\end{proof}
\begin{prop}\label{efru}
Let $x$ has type {\upshape (rm ur)}.
Then $\varepsilon_v(x)=2^{-1}(1-q_v^{-2})(1-q_v^{-1})$.
\end{prop}
\begin{proof}
By Lemma \ref{Dcontainedru}, the set $\cK_vx=G_{\co_v}x$ is equal to
$\#(G_{\co_v}/\prj_1^{-1}(G_{x\,\co_v/\gp_v}))$ number of
disjoint copies of $\cD$.
Since
\begin{equation*}
G_{\co_v}/\prj_1^{-1}(G_{x\,\co_v/\gp_v})\cong
G_{\co_v/\gp_v}/G_{x\,\co_v/\gp_v},
\end{equation*}
by Lemma \ref{stabru} we have
\begin{equation*}
\begin{aligned}
\vol(\cK_vx)
&=\vol(\cD)\cdot\frac{\#(G_{\co_v/\gp_v})}{\#(N_{x\,\co_v/\gp_v})}\\
&=q_v^{-8}\cdot
\frac{(q_v^2(q_v^2-1))^2\cdot(q_v^2-1)(q_v^2-q_v)}{(q_v^2-1)^2}\\
&=(1-q_v^{-2})(1-q_v^{-1}).
\end{aligned}
\end{equation*}
Now the proposition follows from Lemma \ref{evru}.
\end{proof}

Next we consider orbits corresponding to quadratic ramified extensions.
From now on to Proposition \ref{efrr},
we assume $x$ has type (rm rm). Let $x=(1,\varpi)$.
Then $\varpi$ is a prime element of $L_v=k_v(\varpi)\cong k_v(x)$.
Let $\tau$ denote the non-trivial
element of ${\rm Gal}(L_v/k_v)$.
Then $F_x(v_1,v_2)=(v_1+\varpi v_2)(v_1+\varpi^\tau v_2)$
is an Eisenstein polynomial and $(\varpi-\varpi^\tau)^2\in\co_v$
generates the relative discriminant
$\Delta_{k_v(x)/k_v}=\gp_v^{\delta_{x,v}}$.

\begin{lem}\label{evrr}
Let $x$ have type {\upshape (rm rm)}.
Then $\varepsilon_v(x)=\vol(\cK_vx)$.
\end{lem}
\begin{proof}
We can prove this lemma exactly the same as Lemma \ref{evru}.
Only the difference is that we can take the generator $\tau$
of $G_{x\,k_v}/G_{x\,k_v}^\circ$ in Proposition \ref{notidentity}
from $\cK_v$ and hence
$\cK_vG_{x\,k_v}=\cK_vG_{x\,k_v}^\circ$.
\end{proof}

We put $n=\delta_{x,v}+2m_v+1$.
As in the case $x$ has type (ur rm) in Section \ref{sec:fin},
we consider the congruence relation of modulo $\gp_v^n$
to compute $\vol(\cK_vx)$.

\begin{defn}
We define $\cD=\{y\in V_{\co_v}\mid y\equiv x\,(\gp_v^n)\}$.
\end{defn}
\begin{lem}\label{Dcontainedrr}
We have $\cD\subset \cK_vx$.
\end{lem}
\begin{proof}
Let $y\in \cD$.
Then as in the proof of Lemma \ref{Dcontainedur},
we have $y\in G_{k_v}x$.
The rest of argument is similar to that of Lemma \ref{Dcontainedru}
and we shall be brief.
Let $y=gx, g\in G_{k_v}$.
By Lemma \ref{orderstab},
multiplying an element of $G_{x\, k_v}^\circ$
if necessary, we may assume that
\begin{equation*}
|\nr(g_{11})|_v=|\nr(g_{12})|_v=1,\quad \text{and hence}\quad |\det(g_2)|_v=1.
\end{equation*}
Since $F_y(v)\in\sym^2\co_v^2$,
we have
\begin{align*}
F_y(1,0)&=\nr(g_{11})\nr(g_{12})\n_{k_v(x)/k_v}(p+q\varpi)\in\co_v,\\
F_y(0,1)&=\nr(g_{11})\nr(g_{12})\n_{k_v(x)/k_v}(r+s\varpi)\in\co_v.
\end{align*}
Hence both $\n_{k_v(x)/k_v}(p+q\varpi)$ and $\n_{k_v(x)/k_v}(r+s\varpi)$
are elements of $\co_v$ and so $p,q,r,s\in\co_v$.
Since $|\det(g_2)|_v=1$, we conclude $g_2\in\gl(2)_{\co_v}$.
Hence $g\in \cK_v$ and the lemma follows.
\end{proof}
\begin{lem}\label{stabrr}
We have
$[G_{x\,\co_v/\gp_v^n}:N_{x\,\co_v/\gp_v^n}]=2q_v^{\delta_{x,v}}$.
\end{lem}
\begin{proof}
We shall count the number of elements of the right coset space
$N_{x\,\co_v/\gp_v^n}\backslash G_{x\,\co_v/\gp_v^n}$.
Let $g'\in G_{x\, \co_v/\gp_v^n}$.
By Lemma \ref{integralring} and Proposition \ref{stbstr}, the right coset 
$N_{x\, \co_v/\gp_v^n} g'$
contains an element $g=(g_1,g_2)$ with $g_1=(g_{11},g_{12})$
of one of the following forms
\begin{enumerate}
\item[(A)] $g_{11}=1+\alpha\theta,g_{12}=1+\beta\theta$,
\item[(B)] $g_{11}=1+\alpha\theta,g_{12}=\beta+\theta$,
\item[(C)] $g_{11}=\alpha+\theta,g_{12}=1+\beta\theta$,
\item[(D)] $g_{11}=\alpha+\theta,g_{12}=\beta+\theta$,
\end{enumerate}
where $\alpha,\beta\in \co_{L_v}/\gp_v^n$, and also they
are determined by the coset $N_{x\, \co_v/\gp_v^n} g'$ only.
We will count the possibilities for $g$ for each of the above cases.
We choose representatives of $\alpha,\beta$ in $\co_{L_v}$
and use the same notation.

From now on we consider the case $v\notin\gMdy$ and $v\in\gMdy$ separately.
We first consider the case $v\notin\gMdy$.
In this case $\delta_{x,v}=1$ and $n=2$.
Also since $2\in\co_v^\times$, 
by changing $\theta$ and $x=(1,\varpi)$ if necessary,
we may assume that $\theta^2\in\co_v^\times$ and $\varpi^\tau=-\varpi$. 
Let $y=(y_1,y_2)=(g_1,1)x$.

First consider the case of (A).
By computation we have
\begin{equation*}
y_1=1+\alpha\beta^\tau\theta^2+(\alpha+\beta)\theta,
\quad
y_2=\varpi(1-\alpha\beta\theta^2)+\varpi(\beta-\alpha)\theta.
\end{equation*}
Since ${}^t(y_1,y_2)\equiv g_2^{-1}\,{}^t(1,\varpi)\;(\gp_v^2)$,
we have $\prj_2(y_1),\prj_2(y_2)\in\co_{L_v}/\gp_v^2$.
Hence
\begin{equation*}
\alpha+\beta\equiv0\;(\gp_v^2)
\qquad\text{and}\qquad
\varpi(\beta-\alpha)\equiv0\;(\gp_v^2).
\end{equation*}
It is easy to see that
there are $q_v$ possibilities for pairs of $(\alpha,\beta)$
modulo $\gp_v^2$ satisfying the above condition.
On the other hand, for each of these pairs,
we have $y_1\equiv1\;(\gp_v^2)$ and $y_2\equiv\varpi\;(\gp_v^2)$,
and hence $(1,g_2)y\equiv x\;(\gp_v^2)$ if and only if
$g_2\equiv I_2\,(\gp_v^2)$.

Next we consider the case (B).
In this case, we have
\begin{equation*}
y_1=\beta+\alpha\theta^2+(1+\alpha\beta^\tau)\theta,
\quad
y_2=\varpi(\beta-\alpha\theta^2)+\varpi(1-\alpha\beta^\tau)\theta.
\end{equation*}
Again since ${}^t(y_1,y_2)\equiv g_2^{-1}\,{}^t(1,\varpi)$, we have
$\prj_2(y_1),\prj_2(y_2)\in\co_{L_v}/\gp_v^2$.
Hence
\begin{equation*}
1+\alpha\beta^\tau\equiv0\;(\gp_v^2)
\qquad\text{and}\qquad
\varpi(1-\alpha\beta^\tau)\equiv0\;(\gp_v^2),
\end{equation*}
but this is impossible since $2\in\co_v^\times$.
Hence
any right coset of $G_{x\, \co_v/\gp_v^2}$
does not contain elements of the form (B).

The remaining two cases are similar.
We can see that there are no possibilities for $g$ of the form (C) and
$q_v$ possibilities for $g$ of the form (D).
These give the desired description for $v\notin\gMdy$.

We next consider the case $v\in\gMdy$.
In this case we may choose $\theta$ so that
$\theta^2=\theta+c$ for some $c\in\co_v^\times$.
Again we let $y=(y_1,y_2)=(g_1,1)x$.

Let us consider the case $g$ is of the form (A).
By computation we have
\begin{equation*}
y_1=1+c\alpha\beta^\tau
	+(\alpha+\beta+\alpha\beta^\tau)\theta,
\quad
y_2=(\varpi+c\alpha\beta^\tau\varpi^\tau)
	+(\alpha\varpi^\tau+\beta\varpi+\alpha\beta^\tau\varpi^\tau)\theta.
\end{equation*}
Hence as before, we need
\begin{equation*}
\beta+\alpha+\alpha\beta^\tau\equiv0\;(\gp_v^n)
\qquad\text{and}\qquad
\beta\varpi+\varpi^\tau(\alpha+\alpha\beta^\tau)\equiv0\;(\gp_v^n).
\end{equation*}
Under the first equation, the second equation is equivalent to
$\beta(\varpi-\varpi^\tau)\equiv0\;(\gp_v^n)$.
If we write $\beta=\beta_1+\beta_2\varpi$
where $\beta_1,\beta_2\in\co_v/\gp_v^n$,
this equation holds if and only if
\begin{equation*}
\begin{cases}
\beta_1,\beta_2\in\gp_v^{n-\delta_{x,v}/2}/\gp_v^n
			& 2\leq\delta_{x,v}\leq 2m_v,\\
\beta_1\in\gp_v^{n-m_v}/\gp_v^n,\;\;
\beta_2\in\gp_v^{n-m_v-1}/\gp_v^n
			& \delta_{x,v}=2m_v+1,
\end{cases}
\end{equation*}
Hence there are
$q_v^{\delta_{x,v}}$ possibilities for $\beta$.
Also for each of these $\beta$, $1+\beta^\tau$ is invertible
and so $\alpha$ is uniquely determined by the first equation.

Then for each of these pairs $(\alpha,\beta)$,
we have $y_1\equiv1\;(\gp_v^n)$ and $y_2\equiv\varpi\;(\gp_v^n)$,
and hence $(1,g_2)y\equiv x\;(\gp_v^n)$ is equivalent to
$g_2\equiv I_2\,(\gp_v^n)$.
Hence there are $q_v^{\delta_{x,v}}$ choice of 
$g$ of the form (A) in all.

The remaining three cases are similar.
There are no possibilities for $g$ of the form (B) and (C),
and $q_v^{\delta_{x,v}}$ choice of $g$ of the form (D).
We have thus proved the lemma.
\end{proof}

\begin{prop}\label{efrr}
Suppose the standard representative
$x$ has type {\upshape (rm ur)}.
Then
$\varepsilon_v(x)
=	2^{-1}|\Delta_{k_v(x)/k_v}|_v^{-1}(1+q_v^{-1})(1-q_v^{-2})^2$.
\end{prop}
\begin{proof}
By Lemma \ref{Dcontainedrr}, the set $\cK_vx=G_{\co_v}x$ is equal to
$\#(G_{\co_v}/\prj_n^{-1}(G_{x\,\co_v/\gp_v^n}))$ number of
disjoint copies of $\cD$.
Since
\begin{equation*}
G_{\co_v}/\prj_n^{-1}(G_{x\,\co_v/\gp_v^n})\cong
G_{\co_v/\gp_v^n}/G_{x\,\co_v/\gp_v^n},
\end{equation*}
by Lemma \ref{stabrr} we have
\begin{equation*}
\begin{aligned}
\vol(\cK_vx)
&=\vol(\cD)\cdot\frac{\#(G_{\co_v/\gp_v^n})}
{2q_v^{\delta_{x,v}}\cdot\#(N_{x\,\co_v/\gp_v^n})}\\
&=q_v^{-8n}\cdot
\frac{\left\{q_v^{4n}(1-q_v^{-2})\right\}^2
		\cdot q_v^{4n}(1-q_v^{-1})(1-q_v^{-2})}
{2q_v^{\delta_{x,v}}\cdot \left\{q_v^{2n}(1-q_v^{-1})\right\}^2}\\
&=2^{-1}q_v^{-\delta_{x,v}}(1+q_v^{-1})(1-q_v^{-2})^2.
\end{aligned}
\end{equation*}
Now the proposition follows from Lemma \ref{evrr}.
\end{proof}

\section{Computation of the local densities at infinite places}\label{sec:inf}

In this section, we compute $\varepsilon_v(x)$ at infinite places.
We assume $v\in\gMi$ in this section.
For the unramified places, the values were already
computed in \cite{kayuII},
and the remaining case is
for places $v\in\gM_\esB\cap\gMi$.
Note that this case does not occur if $v\in\gMc$ and
that for these places $V_{k_v}^\sst$ is the single $G_{k_v}$-orbit.
In the computation we need to know the $8\times 8$ Jacobian determinant
associated with the map $g\mapsto gx$ in some coordinate system.
This calculation was carried out using
the \textsc{Maple} computer algebra package \cite{maple}.

\begin{prop}
Let $v\in\gM_\esB\cap \gMi$.
Then $\varepsilon_v(x)=\pi^3/2$.
\end{prop}

\begin{proof}

Since $|P(x)|_v=1$ if $x$ is a standard representative,
$\varepsilon_v(x)=b_{x,v}^{-1}$.
We proved in Proposition \ref{bxindep} that if $y\in G_{k_v}x$
then $b_{y,v}=b_{x,v}$.
Therefore we will compute $b_{x,v}$ for $x=(1,\sqrt{-1})$
instead of the standard representative.

We define
\begin{equation*}
\iota\colon \C^\times\longrightarrow \gl(2)_\R
\quad
\text{as}
\quad
t\longmapsto\twtw {\Re(t)}{\Im(t)}{-\Im(t)}{\Re(t)},
\end{equation*}
which is an injective homomorphism.
Then the isomorphism in Proposition \ref{stbstrk}
can be expressed as
\begin{equation*}
\psi_{\sqrt{-1}}\colon
(\C^\times)^2\longrightarrow G_{x\,\R}^\circ,
\qquad
(s_1,s_2)\longmapsto (s_1,s_2,\iota(s_1s_2)^{-1}).
\end{equation*}
Recall that the measure $dg_{x,v}''$ on $G_{x\,\R}^\circ$
was defined as the pushout measure of $\md s_1\md s_2$.

Let $\bbh'=\{(u+j)s\mid u\in\C,s\in\C^\times\}$
and $\cD=\bbh'\times\bbh'\times\gl(2)_{\R}$.
Then $\cD$ is an open dense subset of $G_{k_v}$ and
we will compare the measures on this set.
Any element $g$ of $\cD$ can be written uniquely as
$g=g_{x,v}'g_{x,v}''$ where
\begin{equation*}
g_{x,v}'=(u_1+j,u_2+j,g_3),
\quad
g_{x,v}''=(s_1,s_2,\iota(s_1s_2)^{-1})
\end{equation*}
with
\begin{equation*}
u_1,u_2\in\C,
\quad
g_3=\twtw {a_{11}}{a_{12}}{a_{21}}{a_{22}}\in\gl(2)_\R,
\quad
\text{and}
\quad
s_1,s_2\in\C,
\end{equation*}
and when $g_{x,v}'$ is written in this form,
$\Re(x_i),\Im(x_i)$ and $a_{ij}$ for $i,j=1,2$ may be
regarded as coordinates on $G_\R/G_{x\, \R}^\circ$.
An easy computation shows that
\begin{equation*}
dg_{x,v}'=
\frac{1}{\pi^3}
\cdot \frac{du_1}{(|u_1|_\C+1)^2}
\cdot \frac{du_2}{(|u_2|_\C+1)^2}
\cdot d\mu(g_3)
\end{equation*}
with respect to these coordinates.
Note that we are setting $du_i$ 
twice the Lebesgue measure on $\C$ as usual
and that we defined $d\mu(g_3)$ to be
$da_{11}da_{12}da_{21}da_{22} / |\det(g_3)|_\R^2$.

We consider the Jacobian determinant of the map
\begin{equation*}
G_\R/G_{x\,\R}^\circ\rightarrow G_\R x,
\qquad
g_{x,v}'\mapsto g_{x,v}'x.
\end{equation*}
To do this, we choose there respective $\R$-coordinates.
For $G_\R/G_{x\, \R}^\circ$, 
we regarded 
$\Re(x_i)$, $\Im(x_i)$ and $a_{ij}$ for $i,j=1,2$
as its $\R$-coordinates.
For $G_\R x$, which is an open subset of $V_\R=\bbh\oplus\bbh$,
by expressing elements of $G_\R x$ as
\begin{equation*}
y=(y_{11}+y_{12}j,y_{21}+y_{22}j),
\qquad
y_{ij}\in\C\ (i,j=1,2),
\end{equation*}
we regard $\Re(y_{ij}),\Im(y_{ij})$ for $i,j=1,2$
as $\R$-coordinates of $G_\R x$.
Then with respect to the coordinate systems above,
the Jacobian determinant of the map is found to be
$4(|u_1|_\C+1)^2(|u_2|_\C+1)^2|\det(g_3)|_\R^2$
by using \textsc{Maple} \cite{maple}.
Note that this map is a double cover since
$[G_{x\,\R}:G_{x\,\R}^\circ]=2$.
As $P(g_{x,v}'x)=\chi(g_{x,v}')P(x)$ and
\begin{equation*}
|\chi(g_{x,v}')|_\R=(|u_1|_\C+1)^2(|u_2|_\C+1)^2|\det(g_3)|_\R^2,
\quad
|P(x)|_\R=4,
\end{equation*}
it follows that the pullback measure of $dy/|P(y)|_\R^2$
to $G_{\R}/G_{x\,\R}^\circ$ is
\begin{equation*}
\frac{1}{2}
\cdot \frac{du_1}{(|u_1|_\C+1)^2}
\cdot \frac{du_2}{(|u_2|_\C+1)^2}
\cdot d\mu(g_3).
\end{equation*}
We note that we chose the measure $dy$ on $V_\R$
to be $2^4$ times to that of product of Lebesgue measures
$\prod_{i,j}[d\{\Re(y_{ij})\}d\{\Im(y_{ij})\}]$.
Comparing this measure and $dg_{x,v}'$, we have
$b_{x,v}=2/\pi^3$ and 
hence the proposition follows.
\end{proof}

Combining
\cite[Propositions 5.2, 5.4]{kayuII}
and the above proposition, we obtain the following.

\begin{prop}\label{ei}
Assume $v\in\gMi$.
Let $x\in V_{k_v}^\sst$ be one of the standard representatives.
\begin{enumerate}[{\rm (1)}]
\item If $v\in\gMr$ then
$\varepsilon_v(x)=\pi^3/2$
for any type of the standard representative.
\item If $v\in\gMc$ then
$\varepsilon_v(x)=4\pi^3$.
\end{enumerate}
\end{prop}

All of these finish the necessary preparations from local theory
and we are now ready to go back to the adelic situation.

\section{The mean value theorem}\label{sec:mea}

In this section, we will deduce our mean value theorem
by putting together the results we have obtained before.
We will see in Proposition \ref{ZetaSeries} that
the global zeta function is approximately
the Dirichlet generating series for the sequence $\gC_L^2$
for quadratic extensions $L$ of $k$ which are embeddable into $\esB$.
If it were exactly this generating series,
the Tauberian theorem would allow us
to extract the mean value of the coefficients
from the analytic behavior of this series.
However, our global zeta function contains
an additional factor in each term.
We will surmount this difficulty
by using the technique called the filtering process,
which was originally formulated by Datskovsky-Wright \cite{dawrb}.

Let $x\in V_k^\sst$.
We define measures $dg_x''$ and $d\ti g_x''$ on
$G_{x\A}^\circ$ and $G_{x\A}^\circ/\ti T_\A$
to be $dg_x''=\prod_{v\in\gM}dg_{x,v}''$ and 
$d\ti g_x''=\prod_{v\in\gM}d\ti g_{x,v}''$,
where we defined $dg_{x,v}''$ and $d\ti g_{x,v}''$ in Section \ref{sec:can}.
We choose a left invariant measure on $G_{\bba}/G_{x\,\bba}^{\circ}$.
Since $G_x^\circ$ is isomorphic to $(\gl(1)_{\ti k(x)})^2$
as an algebraic group over $k$,
the first Galois cohomology set $H^1(k',G_x^\circ)$ is trivial
for any field $k'$ containing $k$.
This implies that the set of $k'$-rational point of 
$G_{k'}/G_{x\,k'}^\circ$ coincides with $(G/G_x^\circ)_{k'}$.
Therefore $G_\A/G_{x\,\A}^\circ=(G/G_x^\circ)_\A$.
Hence if we let $dg_x'=\prod_v dg_{x,v}'$
(we defined $dg_{x,v}'$ in Section \ref{sec:loc}),
then this defines
an invariant measure on $G_{\bba}/G_{x\,\bba}^{\circ}$.
We have $dg=dg_x'dg_x''$
since $dg_v=dg_{x,v}'dg_{x,v}''$ for all $v$,
and hence $d\ti g=dg_x'd\ti g_x''$.

We first determine the volume of
$G_{x\,\bba}^{\circ}/{\ti T}_{\bba} G_{x\,k}^{\circ}$ under $d\ti g_x''$,
which is the weighting factor of the Dirichlet series
arising from our global zeta function.

\begin{prop}\label{object}
Suppose {\rm $x\in V^{\sst}_k$}.
Then
the volume of $G^{\circ}_{x\,\A}/\ti T_{\A}G^{\circ}_{x\,k}$ 
with respect to the measure {\rm $d\ti g_x''$}
is $(2\gC_{k(x)}/\gC_{k})^2$.
\end{prop}

\begin{proof}
Identifying $\ti T$ with $(\gl(1)_{k})^2$ and 
$G^{\circ}_x$ with $(\gl(1)_{k(x)})^2$, we define
$\ti T^0_{\A}=(\A^0)^2$ and $G^{\circ0}_{x\,\A}=(\A_{k(x)}^0)^2$.
Let $\md  \ti t^0$ and
$d g_x''{}^0$ be the measures on 
$\ti T^0_{\A}$ and $G^{\circ0}_{x\,\A}$, such that 
$d g_x'' = \md \lam_1\md \lam_2 d g_x''{}^0$, 
$\md  \ti t = \md \lam_1\md \lam_2 \md  \ti t^0$ for  
\begin{equation*}
g_x'' = (\lmb{1}_{k(x)},\lmb{2}_{k(x)}) g_x''{}^0,\quad
\ti t = (\lmb{1}_{k},\lmb{2}_{k}) \ti t^0
.\end{equation*}
Note that if $\lam\in \mr$ then
the absolute value of $\lamb_{k}$ as an idele of 
$ k(x)$ is $\lam^2$.  
Therefore, $d g_x'' = 2^2\md \lam_1\md\lam_2 d g_x''{}^0$
for $g_x'' = (\lmb1_{k},\lmb2_{k})g_x''{}^0$.  
Since $d g_x'' = 
\md  \ti td \ti g_x''$ this implies that
$2^2d g_x''{}^0 = 
\md  \ti t^0d \ti g_x''$. 
Therefore
\begin{equation*}
\begin{aligned}
2^2 \int_{G^{\circ0}_{x\,\A}/G^{\circ}_{x\,k}} dg_x''{}^0
& = \int_{G^{\circ0}_{x\,\A}/G^{\circ}_{x\,k}\ti T_{\A}^0} d\ti g_x'' 
\int_{\ti T^0_{\A}/\ti T_k} \md \ti t^0 \\
& = \vol(G^{\circ0}_{x\,\A}/\ti T_{\A}^0G^{\circ}_{x\,k}) 
\int_{\ti T^0_{\A}/\ti T_k} \md \ti t^0 \,.
\end{aligned}
\end{equation*}
Since 
\begin{equation*}
\int_{G^{\circ0}_{x\,\A}/G^{\circ}_{x\,k}} 
d g_x''{}^0
= \gC_{k(x)}^2\text{\qquad and \qquad}
\int_{\ti T^0_{\A}/\ti T_k} \md  \ti t^0 
= \gC_{k}^2,
\end{equation*}
this proves the proposition.  
\end{proof}

For $x\in V_k^\sst$ and $\Phi=\otimes\Phi_v\in\cS(V_{\A})$ we define the
\emph{orbital zeta function} of $x$ to be
$Z_{x}(\Phi,s)=\prod_{v\in\gM} Z_{x,v}(\Phi_v,s)$.
Note that we defined $Z_{x,v}(\Phi_v,s)$ in Section \ref{sec:loc}.
If $x$ lies in
the orbit of the standard representative $\omega_{v,x}$ in $V_{k_v}^{\sst}$
then we shall write
$\Xi_{x,v}(\Phi_v,s)= Z_{\omega_{v,x},v}(\Phi_v,s)$ and 
$\Xi_{x}(\Phi,s)=\prod_{v\in\gM}\Xi_{x,v}(\Phi_v,s)$. 
We call
$\Xi_x(\Phi,s)$ the 
\emph{standard orbital zeta function}.  
\begin{prop}\label{standardorbital}
For $x\in V_k^\sst$ and $\Phi=\otimes\Phi_v\in\cS(V_{\A})$
we have 
\begin{equation*}
Z_x(\Phi,s)=\esN(\Del_{k(x)/k})^{-s}\Xi_{x}(\Phi,s)\,.
\end{equation*}
\end{prop}
\begin{proof}
By Proposition \ref{sameorbitzeta}, we have
\begin{equation*}
Z_{x}(\Phi,s) =
\left(\prod_{v\in\gM}\frac{|P(\omega_{v,x})|_v}{|P(x)|_v}\right)^s
\Xi_{x}(\Phi,s).
\end{equation*}
Since $P(x)\in k^\times$, we have $\prod_{v\in\gM}|P(x)|_v=1$
by the Artin product formula.
Also since $P(\omega_{v,x})$ generate the local discriminant
$\Delta_{k_v(x)/k_v}$ of $k_v(x)$ if $v\in\gMf$
and $|P(\omega_{v,x})|_v=1$ if $v\in\gMi$,
we have
\begin{equation*}
\prod_{v\in\gM}|P(\omega_{v,x})|_v
=\prod_{v\in\gMf}|P(\omega_{v,x})|_v
=\prod_{v\in\gMf}|\Delta_{k_v(x)/k_v}|_v
=\esN(\Delta_{k(x)/k})^{-1}.
\end{equation*}
Thus we have the proposition.
\end{proof}

\begin{prop}\label{ZetaSeries}
If $\Phi=\otimes\Phi_v\in\cS(V_{\A})$ then we have
\begin{equation*}
Z(\Phi,s)=\frac{2}{\gC_k^2}\sum_{x\in G_k\backslash V_k^\sst}
\frac{\gC_{k(x)}^2}{\esN(\Del_{k(x)/k})^s}\Xi_x(\Phi,s).
\end{equation*}
\end{prop}
\begin{proof}
By the usual modification, we have
\begin{equation*}
Z(\Phi,s)=\sum_{x\in G_k\backslash V_k^\sst}
[G_{x\,k}:G_{x\,k}^\circ]^{-1}
\int_{G^{\circ}_{x\,\A}/\ti T_{\A}G^{\circ}_{x\,k}}d\ti g_x''
\int_{G_\A/G_{x\,\A}^\circ}|\chi(g_x')|^s\Phi(g_x'x)dg_x'.
\end{equation*}
For each $x\in G_k\backslash V_k^\sst$,
the last integral in the above equation is equal to $Z_x(\Phi,s)$
since $\Phi=\otimes_v\Phi_v$ and $dg_x'=\prod_v dg_{x,v}'$.
Now the proposition follows from
$[G_{x\,k}:G_{x\,k}^\circ]=2$ and
Propositions \ref{object}, \ref{standardorbital}.
\end{proof}

We are now ready to describe the filtering process. This process
was originally used in \cite{dawrb} and was also used in
\cite{kayuI}.
Since our situation is quite similar to \cite{kayuI},
we follow this reference.

We fix a finite set $S\supseteq S_0$ of places of $k$.
Let $T$ denote any finite subset $T\supseteq S$ of $\gM$.
Let
$\Xi_{x,v}(s)=\Xi_{x,v}(\Phi_{v,0},s)$.
(We defined $S_0$ and $\Phi_{v,0}$ in Section \ref{sec:loc}.)
\begin{defn}
For each finite subset $T\supseteq S$ of $\gM$, we define
\begin{equation*}
\Xi_{x,T}(s)=\prod_{v\notin T}\Xi_{x,v}(s)
\qquad\text{and}\qquad
L_T(s)=\prod_{v\notin T}L_v(s),
\end{equation*}
where $L_v(s)$ is as in Proposition \ref{estimate}.
\end{defn}
By Proposition \ref{estimate}, we have the following.
\begin{prop}
Both $\Xi_{x,T}(s)$ and $L_T(s)$ are Dirichlet series.
We let
\begin{equation*}
\Xi_{x,T}(s)=\sum_{m=1}^{\infty}\frac{a_{x,T,m}^*}{m^s}
\qquad\text{and}\qquad
L_T(s)=\sum_{m=1}^{\infty}\frac{l_{T,m}^*}{m^{s}}.
\end{equation*}
Then $0\leq a_{x,T,m}^*\leq l_{T,m}^*$ for all $m$ and $a_{x,T,1}^*=1$.
Also $L_T(s)$ converges absolutely and locally uniformly
in the region $\re(s)>3/2$.
\end{prop}

We consider $T$-tuples
$\omega_T=(\omega_v)_{v\in T}$ where each $\omega_v$ is one of the
standard orbital representatives for the orbits in $V_{k_v}^{\sst}$.
If $x\in V_{k}^{\sst}$ and $x\in G_{k_v}\omega_v$
then we write $x\approx \omega_v$ and
if $x\approx\omega_v$ for all $v\in T$
then we write $x\approx \omega_{T}$.

\begin{defn}\label{xisdefn}
We define
\begin{equation*}
\xi_{\omega_T}(s)=\sum_{x\in G_k\backslash V_k^\sst,x\approx\omega_T}
\frac{\gC_{k(x)}^2}{\esN(\Del_{k(x)/k})^s}\Xi_{x,T}(s)
\end{equation*}
and
\begin{equation*}
\xi_{\omega_S,T}(s)=\sum_{x\in G_k\backslash V_k^\sst,x\approx\omega_S}
\frac{\gC_{k(x)}^2}{\esN(\Del_{k(x)/k})^{s}}\Xi_{x,T}(s)\,,
\end{equation*}
which is the sum of $\xi_{\omega_T}(s)$ over all
$\omega_T=(\omega_v)_{v\in T}$ which extend the fixed $S$-tuple
$\omega_S$.
\end{defn}

The following lemma is exactly the same as
\cite[Lemma 6.17]{kayuI} and we omit the proof.
\begin{lem}\label{testfunction}
Let $v\in\gM$, $x\in V_{k_v}^{\text{\upshape ss}}$ and $r\in\C$.
Then there exists a $\cK_v$-invariant
Schwartz-Bruhat function $\Phi_v$ such that the support of
$\Phi_v$ is contained in $G_{k_v}x$, $Z_{x,v}(\Phi_v,s)$ is an
entire function and $Z_{x,v}(\Phi_v,r)\neq0$.
\end{lem}

\begin{prop}\label{xisresidue}
Let $T\supseteq S$ be a finite set of places of $k$ and $\omega_T$
be a $T$-tuple, as above. The Dirichlet series $\xi_{\omega_T}(s)$
has a meromorphic continuation to the region $\re(s)>3/2$. Its only
possible singularity in this region is a simple pole at $s=2$ with
residue
\begin{equation*}
\cR_2
\prod_{v\in\gM_\esB\cap\gMf}(1-q_v^{-1})^2
\prod_{v\in T}\varepsilon_v(\omega_v),
\end{equation*}
where
\begin{equation*}
\cR_2=
\Delta_k^{-5/2}\gC_k^3Z_k(2)^3.
\end{equation*}
\end{prop}
\begin{proof}
For each $v\in T$ we choose $\cK_v$-invariant Schwartz-Bruhat function
$\Phi_v$ such that
$\supp(\Phi_v)\subseteq G_{k_v}\omega_v$. Let
$\Phi=\bigotimes_{v\in T}\Phi_v\otimes\bigotimes_{v\notin
T}\Phi_{v,0}\in\cS(V_{\A})$. For $v\in T$ we have
$\Xi_{x,v}(\Phi_v,s)=0$ unless $x\approx\omega_v$ and hence
by Proposition \ref{ZetaSeries} we have
\begin{align*}
Z(\Phi,s)
&=	\frac{2}{\gC_k^2}
		\sum_{x\in G_k\backslash V_k^\sst,x\approx\omega_T}
			\frac{\gC_{k(x)}^2}{\esN(\Del_{k(x)/k})^s}
			\bigg(\prod_{v\in T}\Xi_{x,v}(\Phi_v,s)\bigg)
			\Xi_{x,T}(s) \\
&=	\frac{2}{\gC_k^2}
		\bigg(\prod_{v\in T}Z_{\omega_v,v}(\Phi_v,s)\bigg)
		\xi_{\omega_T}(s).
\end{align*}
Using Lemma \ref{testfunction}
and Theorem \ref{globalanalyticproperties}, this formula implies
the first statement.
Also by the equation just before Proposition \ref{sameorbitzeta}, we have
\begin{equation*}
Z_{\omega_v,v}(\Phi,2)
=\varepsilon_v(\omega_v)^{-1}\int_{V_{k_v}}\Phi_v(x)dx_v.
\end{equation*}
Since $\int_{V_{k_v}}\Phi_{v,0}(x)\,dx_v=1$ for $v\notin T$,
by Theorem \ref{globalanalyticproperties}
we have the residue of $\xi_{\omega_T}(s)$.
\end{proof}
As a corollary to this proposition, we obtain the following.
\begin{cor}\label{xiSTresidue}
The Dirichlet series $\xi_{\omega_S,T}(s)$ has a meromorphic
continuation to the region $\re(s)>3/2$. Its only possible
singularity in this region is a simple pole at $s=2$ with residue
\begin{equation*}
\cR_2
\prod_{v\in\gM_\esB\cap\gMf}(1-q_v^{-1})^2
\prod_{v\in S}\varepsilon_v(\omega_v)
\cdot
\prod_{v\in T\setminus S}E_v.
\end{equation*}
\end{cor}
We are now ready to prove a preliminary version of the density theorem.
Since the proof is exactly same as that of \cite[Theorem 6.22]{kayuI},
we omit it.
Note that by Proposition \ref{ldsf},
the product
$\prod_{v\in\gM} E_v$ converges to a positive number.
\begin{thm}\label{prelim-main}
Let $S\supset S_0$ be a finite set of places of $k$ and
$\omega_S$ be an $S$-tuple of standard orbital representatives.
Then
\begin{equation*}
\lim_{X\to\infty}\frac{1}{X^2}
\sum_{\substack{x\in G_k\backslash V_k^\sst,x\approx\omega_S \\
\esN(\Del_{k(x)/k})\leq X}}\gC_{k(x)}^2=
\frac12\cR_2
\prod_{v\in\gM_\esB\cap\gMf}(1-q_v^{-1})^2
\prod_{v\in S}\varepsilon_v(\omega_v)\cdot
\prod_{v\notin S} E_v\,.
\end{equation*}
\end{thm}

We will rewrite Theorem
\ref{prelim-main} as a mean value theorem
for the square of class number times regulator
of quadratic extensions.
Let $S\supset\gMi$ be a finite set of places.
As in Section \ref{sec:int},
we let $L_S=(L_v)_{v\in S}$ be an $S$-tuple where
each $L_v$ is a separable quadratic algebra of $k_v$,
and put
\begin{align*}
\esQ(L_S)
&=	\{F\mid \text{$[F:k]=2$,
		$F\otimes k_v\cong L_v$ for all $\ v\in S$}\},\\
\esQ(L_S,X)
&=	\{F\in \esQ(L_S)\mid\esN(\Delta_{F/k})\leq X\}
\end{align*}
where $X$ is a positive number.
%
%
%
To state our main theorem,
we define the constants as follows.

\begin{defn}\label{gii}
\begin{enumerate}[{\rm (1)}]
\item
Let $v\in\gMf$ and $L_v$ a separable quadratic algebra over $k_v$.
We put
\begin{equation*}
\gii_v(L_v)
=
\begin{cases}
	2^{-1}(1+q_v^{-1})(1-q_v^{-2})
		&	L_v\cong k_v\times k_v,\\
	2^{-1}(1-q_v^{-1})^3
		&	L_v\text{ is quadratic unramified},\\
	2^{-1}|\Delta_{L_v/k_v}|_v^{-1}(1-q_v^{-1})(1-q_v^{-2})^2
		&	L_v\text{ is quadratic ramified}.\\
\end{cases}
\end{equation*}
\item
Let $S\supset\gMi$. For a $S$-tuple $L_S=(L_v)_{v\in S}$
of separable quadratic algebras, we define
\begin{align*}
\gii_\infty (L_S)&=2^{-r_1(L_s)}(2\pi)^{-r_2(L_s)}
\end{align*}
by assuming that
\begin{align*}
r_1(L_S)&=\#\{v\in\gMr\mid L_v\cong \R\times\R\}\times 2,\\
r_2(L_S)&=\#\{v\in\gMr\mid L_v\not\cong \R\times\R\}+2r_2.
\end{align*}
\item
For $v\in\gMf$, we put
\begin{equation*}
\gE_v=1-3q_v^{-3}+2q_v^{-4}+q_v^{-5}-q_v^{-6}.
\end{equation*}
Also we define
\begin{equation*}
\cR_k=2^{-(r_1+r_2+1)}e_k^2\gC_k^3.
\end{equation*}
\end{enumerate}
\end{defn}
Note that the constants $\gii_v(L_v)$ are
$(1-q_v^{-2})^{-1}$ times those of we have listed
in Propositions \ref{efuu}, \ref{efur} and
$(1-q_v^{-1})^2(1-q_v^{-2})^{-1}$times those of we have evaluated
in Propositions \ref{efru}, \ref{efrr}.

The following theorem is a main result of this paper.

\begin{thm}\label{maintheorem}
Let $S\supset\gMi$ and $L_S=(L_v)_{v\in S}$ an $S$-tuple.
Assume there are at least 2 places $v$ such that $L_v$ are fields.
Then we have
\begin{equation*}
\lim_{X\to\infty}\frac{1}{X^2}
\sum_{F\in\esQ(L_S,X)}
h_F^2R_F^2
=	\cR_k\Delta_k^{1/2}\zeta_k(2)^2
	\gii_\infty(L_S)^2
	\prod_{v\in S\cap \gMf}\gii_v(L_v)
	\prod_{v\notin S}\gE_v.
\end{equation*}
\end{thm}
\begin{proof}
We choose $v_1,v_2\in S$ so that
$L_{v_1},L_{v_2}$ are fields.
We take the quaternion algebra $\esB$ of $k$ so that
$\gM_\esB=\{v_1,v_2\}$, which is possible by the Hasse principle.
We consider the prehomogeneous vector space $(G,V)$ for this $\esB$.
Since the set of $k_v$-rational orbits
$G_{k_v}\backslash V_{k_v}^\sst$ corresponds
to the set of all quadratic extensions of $k_v$ if $v\in\gM_\esB$ and
to the set of all separable quadratic algebras of $k_v$ if $v\notin\gM_\esB$,
we can take a $S$-tuple $\omega_S=(\omega_v)_{v\in S}$
of standard orbital representatives so that each $\omega_v$
corresponds to $L_v$.
We claim that if a quadratic extension $F$ of $k$ satisfies
$F\in\esQ(L_S)$ then there exists $x\in V_k^\sst$ so that $F\cong k(x)$.
In fact, if $F\in\esQ(L_S)$ then $F\otimes k_{v_i}\cong L_{v_i}$
is embeddable into $\esB_{v_i}$ for $i=1,2$.
Since $\esB_v\cong{\rm M}(2,2)_{k_v}$ for $v\notin\gM_\esB$,
this shows that $F\otimes \esB_v\cong{\rm M}(2,2)_{F\otimes k_v}$
for all $v$ and by the Hasse principle
we have $F\otimes \esB\cong{\rm M}(2,2)_F$.
Hence $F$ is embeddable into $\esB$ and so
by Proposition \ref{rod}, there exists $x\in V_k^\sst$ such that
$F\cong k(x)$.

Therefore, applying Theorem \ref{prelim-main} for $\omega_S$,
we obtain
\begin{equation}\label{gceq}
\lim_{X\to\infty}\frac{1}{X^2}
\sum_{F\in\esQ(L_S,X)}
\gC_{F}^2=
\frac12\cR_2
\prod_{v\in\gM_\esB\cap\gMf}(1-q_v^{-1})^2
\prod_{v\in S}\varepsilon_v(\omega_v)\cdot
\prod_{v\notin S} E_v\,.
\end{equation}
We consider the value $\gC_F^2$.
Let $r_1(F)$ and $r_2(F)$ be the number of set of real places
and complex places, respectively.
Then if $F\in\esQ(L_S)$ we immediately see $r_i(F)=r_i(L_S)$ for $i=1,2$.
Also one can easily see that $e_F=e_k$ all but finitely-many
quadratic extensions $F$ of $k$.
This finite exceptions may be ignored in the limit, and
we have
\begin{equation*}
\gC_F^2=\gii_\infty(L_S)^{-2}e_k^{-2}h_F^2R_F^2.
\end{equation*}
for almost all $F\in\esQ(L_S)$.
Let us consider the right hand side of \eqref{gceq}.
By Proposition \ref{ei} and the definition of $Z_k(s)$, we have
\begin{align*}
\frac12\cR_2\prod_{v\in\gMi}\varepsilon_v(\omega_v)
&=\frac{\gC_k^3}{2\Delta_k^{5/2}}
	\left(\frac{\Delta_k}{\pi^{r_1}(2\pi)^{r_2}}\zeta_k(2)\right)^3
(\frac{\pi^3}{2})^{r_1}(4\pi^3)^{r_2}\\
&=\frac{1}{2^{r_1+r_2+1}}\Delta_k^{1/2}\gC_k^3\zeta_k(2)^3
=e_k^{-2}\cR_k\Delta_k^{1/2}\zeta_k(2)^3.
\end{align*}
Since $\zeta_k(2)=\prod_{v\in\gMf}(1-q_v^{-2})^{-1}$
and $E_v=(1-q_v^{-2})\gE_v$, \eqref{gceq}
turns out that
\begin{multline}\label{hreq}
\lim_{X\to\infty}\frac{1}{X^2}
\sum_{F\in\esQ(L_S,X)}
h_F^2R_F^2
=	\Delta_k^{1/2}\zeta_k(2)^2\cR_k\gii_{\infty}(L_S)\\
\times
\prod_{v\in (S\cap\gMf)\setminus\gM_\esB}
	\frac{\varepsilon_v(\omega_v)}{1-q_v^{-2}}
\prod_{v\in S\cap\gMf\cap\gM_\esB}
	\frac{(1-q_v^{-1})^2\varepsilon_v(\omega_v)}{1-q_v^{-2}}
\prod_{v\notin S}
	\gE_v.
\end{multline}
As in the observation after Definition \ref{gii},
one can see that
$\varepsilon_v(\omega_v)/(1-q_v^{-2})=\gii_v(L_v)$
for $v\in\gMf\setminus\gM_\esB$ and
$(1-q_v^{-1})^2\varepsilon_v(\omega_v)/(1-q_v^{-2})=\gii_v(L_v)$
for $v\in\gMf\cap\gM_\esB$.
Hence we obtain the desired description.
\end{proof}

\begin{rem}
Let $S\supset\gMi$ and $L_S=(L_v)_{v\in S}$
any $S$-tuple of separable quadratic algebras.
For a finite set $T$ of places $L$ of $k$, 
let $\esQ_T$ be the set of quadratic extensions $L$ of $k$
so that $L$ does not split at least two places of $T$.
Then by Theorem \ref{maintheorem}, 
for any $T$ so that $T\cap S=\emptyset$, we could see that
\begin{multline*}
\lim_{X\to\infty}\frac{1}{X^2}
\sum_{F\in\esQ_T,F\in\esQ(L_S,X)}
h_F^2R_F^2\\
=	\cR_k\Delta_k^{1/2}\zeta_k(2)^2
	\gii_\infty(L_S)^2
		\left(
			\prod_{v\in T}\gE_v
		-	\prod_{v\in T}\frac{(1+q_v^{-1})(1-q_v^{-2})}{2}
		\right)			
	\prod_{v\in S\cap \gMf}\gii_v(L_v)
	\prod_{v\notin (S\cup T)}\gE_v
\end{multline*}
and hence
\begin{equation*}
\lim_{T\nearrow(\gM\setminus S)}\lim_{X\to\infty}\frac{1}{X^2}
\sum_{F\in\esQ_T,F\in\esQ(L_S,X)}
h_F^2R_F^2
=	\cR_k\Delta_k^{1/2}\zeta_k(2)^2
	\gii_\infty(L_S)^2
	\prod_{v\in S\cap \gMf}\gii_v(L_v)
	\prod_{v\notin S}\gE_v.
\end{equation*}
If we could change the order of limits
in the right hand side of the above formula,
we can obtain the statement of Theorem \ref{maintheorem}
for unconditional $S$-tuples also.
But to assert the statement is true,
we probably have to know the principal part at the rightmost pole 
of the global zeta function for $\esB={\rm M}(2,2)$,
which is an open problem.
We conclude this section with this conjecture.
\end{rem}
\begin{conj}\label{unconditional}
The statement of Theorem \ref{maintheorem}
also holds for any unconditional $S$-tuple $L_S$.
\end{conj}

\section{The correlation coefficient}
\label{sec:cor}

In this section, we define the correlation
coefficient of class number times regulator
of certain families of quadratic extensions,
and give the value in some cases.
The author would like to thank A. Yukie,
who suggested to consider on this topic.

We fix a quadratic extension $\ti k$ of $k$.
For any quadratic extension $F$ of $k$ other than $\ti k$,
the compositum $F$ and $\ti k$ contains exactly
three quadratic extensions of $k$.
Let $F^\ast$ denote the quadratic extension other than $F$ and $\ti k$.
Note that if we write
$\ti k=k[x]/(x^2-\alpha)$ and $F=k[x]/(x^2-\beta)$
where $\alpha,\beta\in k$ then $F^\ast=k[x]/(x^2-\alpha\beta)$.
As in Section \ref{sec:mea},
let $S$ always denote the finite set of places of $k$ containing $\gMi$ and
$L_S=(L_v)_{v\in S}$ an $S$-tuple of
separable quadratic algebra $L_v$ of $k_v$.
\begin{defn}
We define
\begin{equation*}
\Cor(L_S)
=\lim_{X\to\infty}
	\frac{\sum_{F\in\esQ(L_S,X)}h_FR_Fh_{F^\ast}R_{F^\ast}}
	{\left(\sum_{F\in\esQ(L_S,X)}h_F^2R_F^2\right)^{1/2}
	\left(\sum_{F\in\esQ(L_S,X)}h_{F^\ast}^2 R_{F^\ast}^2\right)^{1/2}}
\end{equation*}
if the limit of the right hand side exists
and call it the \em{correlation coefficient}.
\end{defn}
The asymptotic behavior of the numerator as $X\to\infty$ 
was investigated by \cite{kayuI,kayuII,kayuIII},
while the denominator is considered in this paper.
Hence we could find the correlation coefficients
for certain types of $\ti k$ and $L_S$.
Let $\gM_{\text{rm}}$, $\gM_{\text{in}}$ and $\gM_{\text{sp}}$
be the sets of finite places of $k$ which are
respectively ramified, inert and split on extension to $\ti k$.
Take any $F\in\esQ(L_S)$ to put $L_v^\ast=F^\ast\otimes k_v$
and $L_S^\ast=(L_v^\ast)_{v\in S}$,
which does not depend on the choice of $F$.
In this section we prove the following theorem.
\begin{thm}\label{cor}
Assume $\gMrm\cap\gMdy=\emptyset$ and $S\supset \gMrm$.
Let $L_S=(L_v)_{v\in S}$ is an $S$-tuple of separable quadratic algebras
such that there are at least two places $v$ with $L_v$ are fields.
Further assume that
there are at least two places $v$ with $L_v^\ast$ are fields.
Then we have
\begin{equation*}
\Cor(L_S)
=	\prod_{v\in\gMin\setminus S}
	\left(1-\frac{2q_v^{-2}}{1+q_v^{-1}+q_v^{-2}-2q_v^{-3}+q_v^{-5}}\right).
\end{equation*}
\end{thm}

We first recall from \cite{kayuI}
the asymptotic behavior of
$\sum_{F\in\esQ(L_S,X)}h_FR_Fh_{F^\ast}R_{F^\ast}$
as $X\to\infty$.
We define the constants as follows.
\newcommand{\gff}{\gf}
\begin{defn}\label{gc}
\begin{enumerate}[{\rm (1)}]
\item
Let $v\in\gMf\setminus\gMrm$
and $L_v$ a separable quadratic algebra over $k_v$.
We define $\gff_v(L_v)$ as follows.
\begin{enumerate}[{\rm (i)}]
\item
If $v\in\gMsp$, then we put $\gff_v(L_v)=\gii_v(L_v)$.
\item
If $v\in\gMin$, then we define
\begin{equation*}
\gff_v(L_v)
=
\begin{cases}
	2^{-1}(1-q_v^{-1})(1+q_v^{-2})
		&	L_v\cong k_v\times k_v\text{ or quadratic unramified},\\
	2^{-1}|\Delta_{L_v/k_v}|_v^{-1}(1-q_v^{-1})(1-q_v^{-4})
		&	L_v\text{ is quadratic ramified}.
\end{cases}
\end{equation*}
\item
If $v\in\gMrm\setminus\gMdy$, then we define
\begin{equation*}
\gff_v(L_v)
=
\begin{cases}
	2^{-1}(1-q_v^{-2})
		&	L_v\cong k_v\times k_v,\\
	2^{-1}(1-q_v^{-1})^2
		&	L_v\text{ is quadratic unramified},\\
	2^{-1}q_v^{-2}(1-q_v^{-2})
		&	L_v\cong \ti k_v,\\
	2^{-1}q_v^{-2}(1-q_v^{-1})^2
		&	L_v\text{ is quadratic ramified and }L_v\not\cong\ti k_v.
\end{cases}
\end{equation*}
\end{enumerate}
\item
For an $S$-tuple $L_S=(L_v)_{v\in S}$
we define
$\gff_\infty(L_S)=\gii_\infty (L_S)$.
\item
For $v\in\gMf\setminus \gMrm$, we put
\begin{equation*}
\gF_v=
\begin{cases}
\gE_v	&	v\in\gMsp,\\
(1+q_v^{-2})(1-q_v^{-2}-q_v^{-3}+q_v^{-4})	&	v\in\gMin.
\end{cases}
\end{equation*}
\end{enumerate}
\end{defn}
Then the following is a refinement of
\cite[Corollary 7.17]{kayuI} in case of $\gMdy\cap\gMrm=\emptyset$.
\begin{prop}\label{innerproduct}
Assume $\gMdy\cap\gMrm=\emptyset$ and $S\supset\gMrm$.
Then the limit
\begin{equation*}
\lim_{X\to\infty}\frac{1}{X^2}
\sum_{F\in\esQ(L_S,X)}h_FR_Fh_{F^\ast}R_{F^\ast}
\end{equation*}
exists and the value is equal to
\begin{equation*}
	\cR_k\zeta_{\ti k}(2)\Delta_{\ti k}^{1/2}\Delta_k^{-1/2}
	\gff_\infty(L_S)\gff_\infty(L_S^\ast)
	\prod_{v\in S\cap\gMf}\gff_v(L_v)\prod_{v\notin S}\gF_v.
\end{equation*}
\end{prop}
\begin{proof}
The only new part is that
we determine the constant $\gff_v(L_v)$
for $v\in\gMdy$ and $L_v$ a quadratic ramified extension solely,
whereas in \cite{kayuI}, the sum of $\gff_v(L_v)$
for $L_v$'s with the same relative discriminants were given.
We consider the constants $\gff_v(L_v)$ for these cases.
For $v\in\gMsp$, we could see from \cite{kayuI} that
Proposition \ref{efur} gives not only $\gii_v(L_v)$
but also the value $\gff_v(L_v)$.
Let $v\in\gMin$.
Then a similar argument from Lemma \ref{evur} to Lemma \ref{stabur}
again leads us to the problem to count the number of
the system of congruence equations considered in \cite[Lemma 4.7]{pcbq},
and the result follows.
Since the argument is much the same as the case of $v\in\gMsp$,
we choose not to include the details here.
\end{proof}

We next consider the second term in the denominator.
%
We compare $\Delta_{L_v^\ast/k_v}$ and $\Delta_{L_v/k_v}$.
For $v\in\gMrm$, we put $\sgn(L_v)=-1$
if $L_v$ is a quadratic ramified extension and $\sgn(L_v)=1$ otherwise.
Then in the case $v\notin\gMrm\cap\gMdy$,
the results are described as follows.
\begin{lem}\label{Delta}
We have
$\Delta_{L_v^\ast/k_v}=\gp_v^{\sgn(L_v)}\Delta_{L_v/k_v}$
if $v\in\gMrm\setminus\gMdy$,
while $\Delta_{L_v^\ast/k_v}=\Delta_{L_v/k_v}$
if $v\in\gMsp\cup\gMin$.
\end{lem}
\begin{proof}
For $v\in\gMsp$ these are obviously coincide since $L_v^\ast=L_v$.
If $v\in\gMrm\setminus\gMdy$,
then $L_v$ is quadratic ramified if and only if $L_v^\ast$ is not.
Also $\Delta_{L_v/k_v}$ is $\gp_v$ if $L_v$ quadratic ramified
and is $\co_v$ otherwise, and the result follows.
We consider the case $v\in\gMin$.
If $L_v$ is not quadratic ramified,
then one of $L_v$ and $L_v^\ast$ is the quadratic unramified extension
and the other is $k_v\times k_v$.
Hence their relative discriminants are concurrent.
Assume $L_v$ is quadratic ramified.
If $v\notin\gMdy$ then $L_v$ and $L_v^\ast$ are
the distinct quadratic ramified extensions
of $k_v$ with relative discriminants $\gp_v$,
and therefore the result follows.
If $v\in\gMdy$ then $\Delta_{L_v^\ast/k_v}=\Delta_{L_v/k_v}$
is a corollary of \cite[Proposition 3.10]{kayuIII}.
Thus we obtained the desired description.
\end{proof}
For an $S$-tuple $L_S$, we define
\begin{equation*}
\Delta_{L_S}=\prod_{v\in\gMrm}q_v^{\sgn(L_v)}.
\end{equation*}
\begin{prop}\label{dualmean}
Assume $\gMrm\cap\gMdy=\emptyset$ and $S\supset\gMrm$.
Let $L_S=(L_v)_{v\in S}$ is a $S$-tuple
such that there are at least two places $v$ with $L_v^\ast$ fields.
Then we have
\begin{equation*}
\lim_{X\to\infty}\frac{1}{X^2}
\sum_{F\in\esQ(L_S,X)}
h_{F^\ast}^2R_{F^\ast}^2
=	\Delta_{L_S}^2\cR_k\Delta_k^{1/2}\zeta_k(2)^2
	\gii_\infty(L_S^\ast)^2\prod_{v\in S\cap \gMf}\gii_v(L_v^\ast)
\prod_{v\notin S}\gE_v.
\end{equation*}
\end{prop}
\begin{proof}
By Lemma \ref{Delta}
we have $\esN(\Delta_{F^\ast/k})=\Delta_{L_S}\esN(\Delta_{F/k})$.
Also by definition, $F\in\esQ(L_S)$ if and only if $F^\ast\in\esQ(L_S^\ast)$.
Hence, $F\in\esQ(L_S,X)$ if and only if
$F^\ast\in\esQ(L_S^\ast,\Delta_{L_S}X)$.
Therefore by applying  $L_S^\ast$ to Theorem \ref{maintheorem},
we have the proposition.
\end{proof}
All of these establish the necessary preparations and now
we go back to the proof of Theorem \ref{cor}.
By Theorem \ref{maintheorem}
and Propositions \ref{innerproduct}, \ref{dualmean}, we have
\begin{equation*}
\Cor(L_S)
=	\esN(\Delta_{\ti k/k})^{1/2}\Delta_{L_S}^{-1}
	\frac{\zeta_{\ti k}(2)}{\zeta_k(2)^2}
	\prod_{v\in S\cap\gMf}
		\frac{\gff_v(L_v)}{\left\{\gii_v(L_v)\gii_v(L_v^\ast)\right\}^{1/2}}
	\prod_{v\not\in S}\frac{\gC_v}{\gE_v}.
\end{equation*}
Note that we used the relation
$\esN(\Delta_{\ti k/k})=\Delta_{\ti k}/\Delta_k^2$.
We naturally regard the right hand side
of the equation above as the Euler product of finite places
$\prod_{v\in\gMf}\alpha_v$.
Then we immediately see $\alpha_v=1$ for $v\in\gMsp$
and
\begin{equation*}
\alpha_v
=	\frac{(1-q_v^{-2})^2}{1-q_v^{-4}}\cdot
	\frac{\gC_v}{\gE_v}
=1-\frac{2q_v^{-2}}{1+q_v^{-1}+q_v^{-2}-2q_v^{-3}+q_v^{-5}}
\end{equation*}
for $v\in\gMin\setminus S$.
Now the remaining task is to verify $\alpha_v=1$ for $v\in S\setminus\gMsp$
and this could be easily carried out by one by one calculation.
For example, if $v\in S\cap\gMin$ and $L_v$ is quadratic ramified,
then we have
\begin{equation*}
\alpha_v
=	\frac{(1-q_v^{-2})^2}{1-q_v^{-4}}\cdot
	\frac	{2^{-1}|\Delta_{L_v/k_v}|_v^{-1}(1-q_v^{-1})(1-q_v^{-4})}
			{2^{-1}|\Delta_{L_v/k_v}|_v^{-1}(1-q_v^{-1})(1-q_v^{-2})^2}
=	1,			
\end{equation*}
and if $v\in S\cap\gMrm$ 
(note that by the assumption $\gMrm\cap\gMdy=\emptyset$,
we have $v\not\in\gMdy$ in this case)
and $L_v$ is quadratic unramified,
then we have
\begin{equation*}
\alpha_v
=	q_v^{(1/2)-1}(1-q_v^{-2})\cdot
	\frac	{2^{-1}(1-q_v^{-1})^2}
			{2^{-1}q_v^{-1/2}(1-q_v^{-1})^2(1-q_v^{-2})}
=	1.
\end{equation*}
The other cases are similar and we omit the routine figuring here.

\begin{rem}
The conditions on $S$-tuple $L_S$ in Theorem \ref{cor}
is to use Theorem \ref{maintheorem}
for $L_S$ and $L_S^\ast$.
If Conjecture \ref{unconditional} is true,
then we could obtain Theorem \ref{cor} for unconditional $L_S$ also.
\end{rem}



\end{document}